\documentclass[10pt]{article}
\usepackage[latin1]{inputenc}
\usepackage{epsfig}
\usepackage{color}
\usepackage[british,english]{babel}
\usepackage{amsthm}
\usepackage{amsmath}
\usepackage{amsfonts}
\usepackage{amssymb}
\usepackage{graphicx}
\newtheorem{corollary}{Corollary}[section]
\newtheorem{definition}[corollary]{Definition}

\newtheorem{lemma}[corollary]{Lemma}
\newtheorem{proposition}[corollary]{Proposition}
\newtheorem{remark}[corollary]{Remark}
\newtheorem{theorem}[corollary]{Theorem}
\newfont{\sBlackboard}{msbm10 scaled 900}

\newcommand{\mylabel}[1]{\label{#1}
            \ifx\undefined\stillediting
            \else \fbox{$#1$}\fi }
\newcommand{\BE}{\begin{equation}}

\newcommand{\EEQ}{\end{equation}}
\newcommand{\rfb}[1]{\mbox{\rm
   (\ref{#1})}\ifx\undefined\stillediting\else:\fbox{$#1$}\fi}

\newfont{\Blackboard}{msbm10 scaled 1200}

\newfont{\roma}{cmr10 scaled 1200}

\def\CC{\rm \hbox{C\kern-.56em\raise.4ex
         \hbox{$\scriptscriptstyle |$}\kern+0.5 em }}

\def \noame{\noalign{\medskip}}
\newcommand{\ep}   {\epsilon}
\usepackage{fancyhdr}

\makeatletter
\def\section{\@startsection {section}{1}{\z@}{-3.5ex plus -1ex minus
    -.2ex}{2.3ex plus .2ex}{\large\bf}}
\makeatother
\def\be{\begin{equation}}
\def\ee{\end{equation}}

\date{ }
\begin{document}
\thispagestyle{empty}
\title{\Large \bf  Modeling of a micropolar thin film flow with rapidly varying thickness and non-standard boundary conditions}\maketitle
\vspace{-2cm}
\begin{center}
Mar\'ia ANGUIANO\footnote{Departamento de An\'alisis Matem\'atico. Facultad de Matem\'aticas. Universidad de Sevilla. 41012-Sevilla (Spain) anguiano@us.es} and Francisco Javier SU\'AREZ-GRAU\footnote{Departamento de Ecuaciones Diferenciales y An\'alisis Num\'erico. Facultad de Matem\'aticas. Universidad de Sevilla. 41012-Sevilla (Spain) fjsgrau@us.es}
 \end{center}

 \renewcommand{\abstractname} {\bf Abstract}
\begin{abstract}
In this paper, we study the asymptotic behavior of the micropolar fluid flow through a thin domain assuming zero Dirichlet boundary condition on the top boundary, which is rapidly oscillating, and  non-standard boundary conditions on the flat bottom. Assuming ``Reynolds roughness regime", in which the thickness  of the domain is very small compared to the wavelenth of the roughness (i.e. a very slight roughness), we rigorously derive  a   generalized Reynolds equation for pressure clearly showing the roughness-induced effects. Moreover, we give expressions for the average velocity and microrotation. 
\end{abstract}
\bigskip\noindent

\noindent {\small \bf AMS classification numbers:}  35B27, 35Q35, 76A05, 76M50, 76A20.  \\

\ \\
\section {Introduction}\label{S1}

Eringen \cite{eringen} in 60's proposed the model of micropolar fluid which has been extensively studied both in the engineering and mathematical literature, due to its practical importance. This model takes into consideration the microstructure of the fluid particles and captures the effects of its rotation, so the micropolar fluid model describes the motion of a lot of real fluids in a better way than the classical (Newtonian) model. Some examples are liquid crystals, animal blood, muddy fluids, certain polymeric fluids or even water in models with small scales. The rotation of the fluid particles is mathematically described by introducing the microrotation field, together with the standard velocity and pressure fields. Accordingly, a new governing equation coming from the conservation of angular momentum.   Thus, the stationary and incompressible micropolar fluid flow is described by the following equations  expressing the balance of momentum, mass and angular momentum (see  Lukaszewicz \cite{Luka})
\begin{equation}\label{Micro_intro2_dimension}
\left\{\begin{array}{rl}
\displaystyle
-(\nu+\nu_r)\Delta {\bf  \overline u} +\nabla \overline p=2\nu_r{\rm rot}({\bf\overline  w} ),\\
\noame 
\displaystyle  {\rm div}({\bf \overline u})=0,\\
\noame
\displaystyle -(c_a+c_d)\Delta{\bf \overline w} +4\nu_r{\bf\overline w} =2\nu_r{\rm rot}({\bf \overline w}).
\end{array}\right.
\end{equation}
We observe that, in view of the application we want to model, we can assume a small Reynolds number and neglect the inertial terms in the governing equations. In the above system, ${\bf \overline u}=({\bf \overline u}'(\bar x), \overline u_{3}(\bar x))$, with ${\bf \overline u}'=(\overline u_1, \overline u_2)$, denotes the velocity vector field, $ \overline p= \overline p(\bar x)$ the scalar pressure, and ${\bf  \overline w}=({\bf \overline w}'(\bar x),  \overline w_{3}(\bar x))$, with ${\bf \overline w}'=(\overline w_1, \overline w_2)$,  the microrotation field.  $\nu$ is the Newtonian viscosity and  $\nu_r, c_a, c_d$ are microrotation viscosities that characterize the
isotropic properties of the fluid (see e.g. \cite{eringen, Luka}).

An important problem in lubrication theory is to derive lower-dimensional models to describe the flow of fluid in a narrow space between two surfaces in relative motion.  To describe the classical setting, we can define a thin  domain   by 
\begin{equation}\label{ThinDomain0}
 \overline \Omega^\ep=\left\{\overline x=(\overline x',\overline x_3)\in \mathbb{R}^3\,:\, \overline x'=(\overline x_1, \overline x_2)\in L \omega, \quad 0<\overline x_3< c h\left({\overline x'\over L}\right)
 \right\},
\end{equation}
with $\omega\subset \mathbb{R}^2$ is an open subset with smooth boundary,   $L$ is  the characteristic length of the domain $\omega$, $c$ stands for the maximum distance between the surfaces, and $h$ represents the thickness of the domain (which is stationary), with $h$ a smooth $\mathbb{R}$-valued function defined on $\mathbb{R}^2$. We consider the ratio $\ep =c/L$ 
as the small parameter.

It has been observed (see e.g. Bayada {\it et al.}  \cite{Bayada_NewModel} and Bayada and Lukaszewicz \cite{BayadaLuc}) that the magnitude of the viscosity coefficients appearing in (\ref{Micro_intro2_dimension}) may influence the effective flow. Thus, it is reasonable to work with the system written in a dimensionless form. In view of that,  we introduce the characteristic velocity $V_0$ of
the fluid, and define
\begin{equation}\label{adimensionalization_change}
\begin{array}{c}
\displaystyle x={\overline x\over L},\quad 
{\bf u}={{\bf \overline u}\over V_0},\quad p={L\over  V_0(\nu+\nu_r)}\overline p,\quad {\bf w}={L\over V_0}{\bf \overline w},\\
\noame
\displaystyle
N^2={\nu_r\over \nu+\nu_r},\quad R_M={c_a+c_d\over \nu+\nu_r}{1\over L^2}.
\end{array}
\end{equation}
 Dimensionless (non-Newtonian parameter $N^2$ characterizes the coupling between the equations for the velocity and microrotation and it is of order $\mathcal{O}(1)$ with respect to $\ep$. The second dimensionless parameter, denoted by $R_M$ is, in fact, related to the characteristic length of the microrotation effects and will be compared with $\ep$. 
 
 In view of the above change, the stationary and linearized micropolar equations in dimensionless form is given by 

\begin{equation}\label{Micro_intro2}
\left\{\begin{array}{rl}
\displaystyle
-\Delta {\bf u} +\nabla p=2N^2{\rm rot}({\bf w} ),\\
\noame 
\displaystyle  {\rm div}({\bf u})=0,\\
\noame
\displaystyle -R_M\Delta{\bf w} +4N^2{\bf w} =2N^2{\rm rot}({\bf w}),
\end{array}\right.
\end{equation}
where ${\bf u}$ denotes the velocity vector field, $p$ the pressure and ${\bf w}$  the microrotation field.  Moreover, the rescaled domain is given by
\begin{equation}\label{ThinDomain}
 \widehat \Omega^\ep=\omega\times (0, \ep h(x')),\quad  x'=(x_1,x_2)\in\omega\subset\mathbb{R}^2,\quad 0<\ep\ll 1,
\end{equation}
where  $\ep h$ represents the thickness of the domain, with $h$  a smooth function and the small parameter $\ep$  is devoted to tend to zero.   To close up the governing problem, one should specify the reasonable boundary conditions for the velocity and microrotation. Thus, prescribing velocity on the bottom $\widehat \Gamma_0=\omega\times \{0\}$ and zero on the top $\widehat \Gamma_1^\ep=\omega\times \{\ep h(x')\}$ (we don't care about the lateral boundary conditions here), i.e.
\begin{equation}
\begin{array}{l}\displaystyle
{\bf u}={\bf s}=({\bf s}',0)=(s_1, s_2, 0)\quad\hbox{on }\widehat\Gamma_0,\qquad
  {\bf u}=0\quad\hbox{on }  \widehat \Gamma_1^\ep,
\end{array}
\end{equation}
 and  imposing that the mean value of the pressure is zero, then Bayada \& Lukaszewicz \cite{BayadaLuc}  gave a mathematical proof of the transition from the micropolar equations (\ref{Micro_intro2}) to a 2D Reynolds-like equation   by using homogenization techiniques when $\ep\to 0$. This equation is  given by
 \begin{equation}\label{Reynolds_classiq_micr}
{\rm div}_{x'}\left(-{h^3(x')\over 1-N^2}\Phi(h(x'),N)\nabla_{x'}  p(x')+  {h(x')\over 2}{\bf s}'\right)=0\quad\hbox{in}\quad \omega,
\end{equation}
where $\Phi(h(x'),N)={1\over 12}+{1\over 4h^2(x')(1-N^2)}-{1\over 4h(x')}\sqrt{{N^2\over 1-N^2}}\coth\left(Nh(x')\sqrt{1-N^2}\right)$.

Throughout the mathematical literature, micropolar fluid models have been extensively studied in recent years by deriving different asymptotic Reynolds-like models taking into account that a surface presents micro-roughness, depending on the ratio between the size of the roughness and the thickness of the domain and on the boundary conditions considered, see for instance  Anguiano \& Su\'arez-Grau \cite{Anguiano_SG_magneto}, Bayada {\it et al.} \cite{Bayada_Gamouana, Bayada_Gamouana2}, Boukrouche \cite{Boukrouche1}, Boukrouche \& Paoli \cite{Boukrouche2}, Boukrouche {\it et al.} \cite{Boukrouche3, Boukrouche4}, Dupuy {\it et al.} \cite{Dupuy1, Dupuy2}, Mahabaleshwar {\it et al.} \cite{Mahabaleshwar},  Marusic-Paloka \cite{Marusic-Paloka, Marusic-Paloka2, Marusic-Paloka3},  Pa\v zanin \cite{Pazanin}, Pa\v zanin \& Radulovi\'c \cite{PazaninRadulovic}, Pa${\rm \check{z}}$anin \& Su\'arez-Grau \cite{Pazanin_SG},  Su\'arez-Grau \cite{SG_porous} and references therein.
\\

One important problem in the lubrication framework concerns to the case of a thin domain with rapidly oscillating boundary given (after adimensionalization) by 
\begin{equation}\label{ROB}
\Omega^\ep=\omega\times (0, h_\ep(x'))\quad \hbox{with}\quad h_\epsilon(x')=\epsilon h(x'/\epsilon^\ell)\quad \forall\,x'\in\omega\subset\mathbb{R}^2,
\end{equation}
where $h_\ep$  represents the thickness of the domain, which is rapidly oscillating. The thickness is given by the small parameter $0<\ep\ll 1$, and the top boundary is rough with small roughness of wavelength described by $0<\ep^\ell\ll 1$, with $\ell\in (0, +\infty)$. Thus, depending on the limit of $\lambda=\lim_{\ep\to 0}\ep^{1-\ell}$,  there exist three characteristic regimes:
\begin{itemize}
\item[--] Reynolds roughness ($\lambda=0$), where the thickness  of the domain is smaller than the wavelenth of the roughness, i.e. a very slight roughness.
\item[--] Stokes roughness ($\lambda=1$), where  the thickness  of the domain and the wavelenth of the roughness are proportional. 
\item[--] High-frequency roughness ($\lambda=+\infty$), where the thickness of the domain is greater than the wavelenth of the roughness, i.e. a very strong roughness.
\end{itemize}
These characteristic regimes have been studied in  Bayada and Chambat \cite{Bayada_Chambat_1988, Bayada_Chambat}  Benhaboucha {\it et. al.} \cite{Benhaboucha}, Fabricius {\it et al.} \cite{Fab2} and Mikeli\'c \cite{Mikelic2}, Pa\v zanin \& Su\'arez-Grau in \cite{Pazanin_SG_two_osci} for  Newtonian fluids, in Pa\v zanin \& Su\'arez-Grau \cite{Pazanin_SG_thermo} and Su\'arez-Grau \cite{SG1} for micropolar fluids, and in  Anguiano \& Su\'arez-Grau \cite{Anguiano_SG} and Nakasato \& Pa${\rm\check{z}}$anin \cite{Nakasato} for generalized Newtonian with power law viscosity law. 
 \\

Note that previous results consider Dirichlet boundary condition for  microrotation. However, Bayada {\it et al.} \cite{Bayada_NewModel, Bayada_NewModel2}  introduced a new boundary condition more general (and physically justified)  at the fluid-solid interface.   This boundary condition, called non-standard or non-zero boundary condition in the literature, is based on the concept of the so-called {\it boundary viscosity} and was originally proposed in Bessonov \cite{Bessonov1, Bessonov2}, and  links the value of the micorotation with the rotation of the velocity in the following manner
\begin{equation}\label{NSBC1}{\bf w}\times {\bf n}={\alpha\over 2}{\rm rot}({\bf u})\times {\bf n},\quad {\bf w}\cdot {\bf n}=0,
\end{equation}
where ${\bf n}$ is the normal unit vector to the boundary. The coefficient  $\alpha$ is defined as a microrotation retardation at the boundary and it is connected with the different viscosity coefficients.  It has been shown experimentally in Jacobson \cite{Jacobson} and L\'eger {\it et al.} \cite{Leger} that there are chemical interactions beween a solid surface and the nearest fluid layer. This had to be taken into account, especially for a non-Newtonian fluid and a very thin-film thickness. This can be done by introducing a viscosity $\nu_b$ near the surface which is different from $\nu$ and $\nu_r$. In \cite{Bessonov2}, it was proposed to define $\alpha$ by means of this boundary viscosity $\nu_b$ by
$$\alpha={\nu+\nu_r-\nu_b\over \nu_r}.$$
This, following \cite{Bessonov2}, it is possible to give phisical limist to $\nu_b$, inducing limits on $\alpha$:
$$0\leq \nu_b\leq \nu+\nu_r\quad\Rightarrow \quad 0\leq \alpha\leq {\nu+\nu_r\over \nu_r}={1\over N^2}.$$
The condition $\alpha=0$ is equivalent to strong adhesion of the fluid particles to the boundary surface, so that they do not rotate relative to the boundary, i.e. ${\bf w}=0$. Thus, from now on, we consider $\alpha>0$, so that the stress tensor and the microrotation are coupled on the boundary. 

Imposing that the mean value of the pressure is zero, the authors in \cite{Bayada_NewModel, Bayada_NewModel2} managed to prove the existence and uniqueness of the corresponding weak solution providing the well-posedness of the governing problem. It is shown that, in such setting, classical no-slip condition for the velocity should be replaced by the the condition allowing the slippage in the tangential direction and retains a non-penetration condition in the normal direction on the boundary as follows
\begin{equation}\label{NSBC}
\beta({\bf u}-{\bf s})\times {\bf n}={R_M\over 2N^2}{\rm rot}({\bf w})\times {\bf n},\quad {\bf u}\cdot {\bf n}=0.
\end{equation}
The value ${\bf u}-{\bf s}$, if not zero,  is a characteristic of a slippage on the boundary, and $\beta$  is an additional parameter which will enable the influence of this new condition to be controlled. Recently, this type of boundary conditions is attracting a lot of attention, as can be seen in various studies, see for instance  Bene\v s  {\it et al.} \cite{Benes}, Bonnivard {\it et al.} \cite{Bonn_Paz_SG2, Bonn_Paz_SG}, Pa\v zanin \cite{PazaninFilomat} and Su\'arez-Grau \cite{SG3}. \\

Coming back  to the classic setting in lubrication, starting from the micropolar equations (\ref{Micro_intro2}) in the thin domain $\widehat\Omega^\ep$ with non-standard boundary conditions on the bottom $\widehat \Gamma_0$ and zero boundary conditions on the top $\widehat \Gamma_1^\ep$,  after the homogenization process, a generalized Reynolds equation  was derived in  Bayada {\it et al.} \cite{Bayada_NewModel, Bayada_NewModel2} of the form 
\begin{equation}\label{Reynolds_classiq_micr_NSBC}
-{\rm div}_{x'}\Big(\Theta_1(x')\nabla_{x'}p(x')+\Theta_2(x'){\bf s}'\Big)=0\quad \hbox{in }\omega,
\end{equation}
 for certain explicit functions $\Theta_1$ and $\Theta_2$, see \cite[Theorem 3.3]{Bayada_NewModel} and  \cite[Section 2]{Bayada_NewModel2}.   Moreover, explicit expressions for the limit velocity and microrotation were given.

The purpose of this paper is to generalize the results in Bayada {\it et al.} \cite{Bayada_NewModel, Bayada_NewModel2}  to the case of rapidly oscillating domain $\Omega^\ep$ given by (\ref{ROB}) and derive a generalized lower-dimensional model influenced by the non-standard boundary conditions and the roughness of the top boundary $\Gamma_1^\ep=\omega\times \{\ep h(x'/\ep^\ell)\}$. We resctrict ourselves to the case of ``Reynolds roughness", i.e.  the case of a rapidly oscillating boundary described by (\ref{ROB}) with $\lim_{\ep\to 0}\ep^{1-\ell}=0$, which will allows us to derive an explicit equation. Accordingly, we consider the micropolar equations (\ref{Micro_intro2}) assuming non-standard boundary conditions (\ref{NSBC1})--(\ref{NSBC}) on the bottom $\Gamma_0$ and zero boundary conditions on the top boundary $\Gamma_1$. As far as the authors know, the flow of a micropolar fluids with non-standard boundary conditions has not been yet considered in the above described setting.  By applying reduction of dimension techniques together with an adaptation of the unfolding method to capture the micro-geometry of the roughness (see Section \ref{sec:unfolding}),  we rigorously derive  a generalized micropolar Reynolds equation given by
$$-{\rm div}_{x'}\left(K\nabla_{x'} p(x')+ L{\bf s}'\right)=0\quad\hbox{in }\omega,
$$
with flow factors $K\in \mathbb{R}^{2\times 2}$ and $L\in\mathbb{R}$ respectively given by
$$K_{ij}=\int_{Z'}\Theta_1(z')\left(\partial_{z_i}q^j(z')+\delta_{ij}\right)dz',\quad i,j=1, 2,\qquad  L=\int_{Z'}\Theta_2(z')\,dz',$$
where the function $q_i$, $i=1,2$ satisfies the following local periodic problem
$$-{\rm div}_{z'}\left(\Theta_1(z')(\nabla_{z'}q_i(z')+{\bf e}_i)+\Theta_2(z'){\bf s}'\cdot {\bf e}_i\right)=0\quad\hbox{in }Z',\quad i,j=1,2,
$$
with $Z'$ is the periodic cell in $\mathbb{R}^2$ and functions $\Theta_1$ and $\Theta_2$ are the same functions as in (\ref{Reynolds_classiq_micr_NSBC}), which are defined later in Lemma 6.3. In addition, we give the expressions for the average of the velocity and microrotation (see Theorem \ref{thm_main}). Since the obtained findings are amenable for the numerical simulations, we believe that it could prove useful in the engineering practice as well.

The paper is structured as follows.  In Section \ref{sec:setting}, we introduce the position of the problem and the existence and uniqueness results. To pass to the limit, it is necessary a rescaling of the problem, which is introduced in Section \ref{sec:rescaling}. In Section \ref{sec:estimates}, we get the {\it a priori} estimates for velocity, microrotation and pressure, extended to a  $\ep$-independent domain. To capture the micro-geometry of the roughness, we introduce in Section \ref{sec:unfolding} the adaptation of the unfolding method, give some properties and derive the limit problem. Finally, we deduce the generalize Reynolds equation in Section \ref{sec:reynolds}. We finish the paper with a list of references. 

\section{Position of the problem}\label{sec:setting}
In this section, we first define the rough thin domain, the necessary sets to study the asymptotic behavior of the solutions and some notation.  Next, we introduce the micropolar equations and the boundary conditions in the thin domain. Finally, we give the variational formulation and the condition for existence and uniqueness  of solution of the problem according to this setting.

\subsection{Definition of the domain and some notation.}  
\begin{figure}[h!]
\begin{center}
\includegraphics[width=15cm]{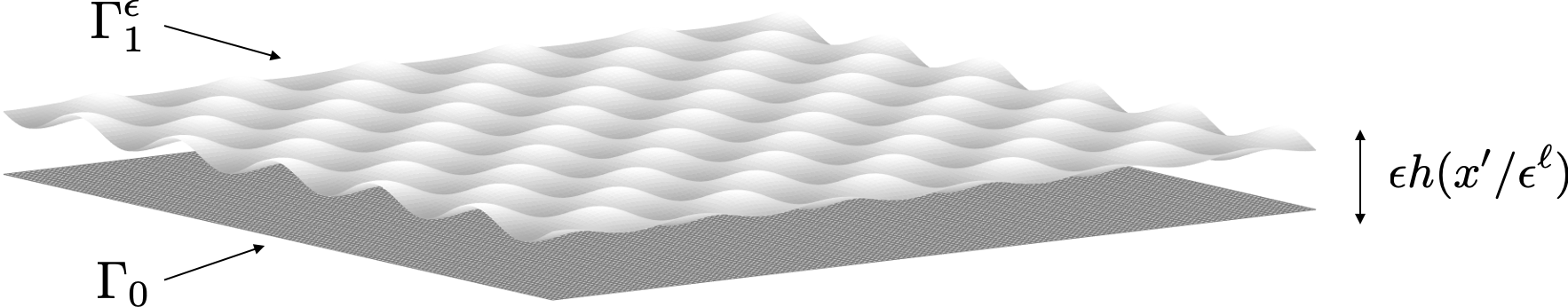}
\caption{Domain $\Omega^\ep$, bottom flat boundary $\Gamma_0$ and top oscillating boundary $\Gamma_1^\ep$}
\label{fig:domain}
\end{center}
\end{figure}

\begin{figure}[h!]
\begin{center}
\includegraphics[width=10cm]{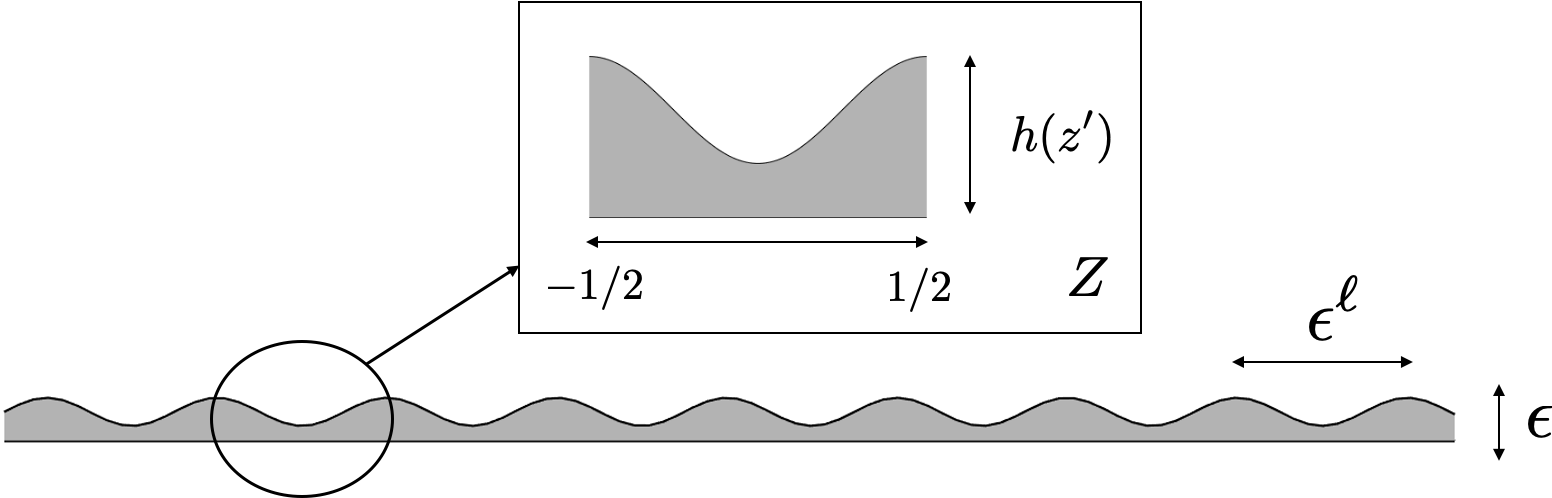}
\caption{Domain $\Omega^\ep$ in 2D and the reference cell Z in 2D}
\label{fig:domain2D}
\end{center}
\end{figure}

We consider $\omega$ a smooth and connected subset of $\mathbb{R}^2$. We define the thin domain with a rapidly varying thickness  by
\begin{equation}\label{Omegaep}
\Omega^\epsilon=\{x=(x',x_3)\in\mathbb{R}^2\times \mathbb{R}\,:\, x'\in \omega,\ 0<x_3< h_\epsilon(x')\},
\end{equation}
where the function $h_\ep(x')= \epsilon h\left(x'/\epsilon^\ell\right)$ represents the real distance between the two surfaces, see Figure \ref{fig:domain}. The small parameter $0<\epsilon\ll 1$ is related to the film thickness and  the small parameter $0<\epsilon^\ell\ll 1$ is the wavelength of the roughness, see Figure \ref{fig:domain2D}. As said in the introduction, we restrict ourselves to the ``Reynolds roughness regime", so we consider that $\epsilon^\ell$ is of order greater or equal than $\epsilon$, i.e. $0<\ell<1$, which implies
\begin{equation}\label{RelationAlpha}
\lim_{\ep\to 0}\ep^{1-\ell}=0.
\end{equation}

\noindent Function $h$ is a positive and smooth function defined for $z'$ in $\mathbb{R}^2$, $Z'$-periodic with $Z'=(-1/2,1/2)^2$ the cell of periodicity in $\mathbb{R}^2$, and there exist $h_{\rm min}$ and $h_{\rm max}$ such that
$$0<h_{\rm min}=\min_{z'\in Z'} h(z'),\quad  h_{\rm max}=\max_{z'\in Z'}h(z')\,.$$

\noindent
We define the boundaries of $\Omega^\epsilon$ as follows
$$\begin{array}{c}
\displaystyle \Gamma_0 =\omega\times\{0\},\quad \Gamma_1^\epsilon=\left\{(x',x_3)\in\mathbb{R}^2\times\mathbb{R}\,:\, x'\in \omega,\ x_3=  h_\epsilon(x')\right\},\quad
\displaystyle\Gamma_{\rm lat}^\epsilon=\partial\Omega_\ep\setminus (\Gamma_0\cup\Gamma_1^\ep).
\end{array}$$

To define the microstructure of the periodicity of the boundary,  we assume that the domain $\omega$ is divided by a mesh of size $\ep^\ell$: for $k'\in\mathbb{Z}^2$, each cell $Z'_{k',\ep^\ell}=\ep^\ell k'+\ep^\ell Z'$. We define $T_\ep=\{k'\in\mathbb{Z}\,:\, Z'_{k',\ep^\ell}\cap\omega\neq \emptyset\}$. In this setting, there exists an exact finite number of periodic sets $Z'_{k',\ep^\ell}$ such that $k'\in T_\ep$.   Also, we define $Z_{k',\ep^\ell}=Z'_{k',\ep^\ell}\times (0,h(z'))$ and  $Z=Z'\times (0,h(z'))$, which is the reference cell in $\mathbb{R}^3$, see Figure \ref{fig:domain2D}. We define the boundaries $\widehat \Gamma_0=Z'\times \{0\}$, $\widehat \Gamma_1=Z'\times\{h(z')\}$. 

\newpage
\noindent Applying a dilatation in the vertical variable, i.e. $z_3=x_3/\ep$, we define the following rescaled sets 
 \begin{equation}\label{domains_tilde}\begin{array}{c}
 \displaystyle \widetilde \Omega^\epsilon=\{(x',z_3)\in\mathbb{R}^2\times \mathbb{R}\,:\, x'\in \omega,\ 0<z_3< h(x'/\ep^\ell)\},\\
 \noame
  \widetilde\Gamma_1^\epsilon=\{(x',z_3)\in\mathbb{R}^2\times \mathbb{R}\,:\, x'\in \omega,\ z_3=h(x'/\ep^\ell)\}\quad \widetilde \Gamma_{\rm lat}^\epsilon=\partial\widetilde \Omega_\ep\setminus (\widetilde \Gamma_0\cup\widetilde \Gamma_1^\ep).
  \end{array}
  \end{equation}
The quantity  $h_{\rm max}$ allows us to define:
\begin{itemize}
\item[--] The extended sets $Q_\ep=\omega\times (0,\ep h_{\rm max})$,  $\Omega=\omega\times  (0, h_{\rm max})$  and $\Gamma_1=\omega\times \{h_{\rm max}\}$.
\item[--]  The extended cube  $\widetilde Q_{k',\varepsilon}=Z'_{k',\varepsilon}\times (0, h_{\rm max})$ for $k'\in\mathbb{Z}^2$.

\item[--]  The extended and restricted basic cell $\Pi=Z'\times (0,h_{\rm max})$.
\end{itemize}

\noindent We denote by $C$ a generic constant which can change from line to line. Moreover, $O_\ep$ denotes a  generic quantity, which can change from line to line, devoted to tend to zero when $\ep\to 0$.\\

\noindent We use the following notation for the partial differential operators
 $$\begin{array}{c}
 \displaystyle \Delta{\bf \varphi}=\Delta_{x'}{\bf \varphi} +\partial_{x_3}^2 {\bf \varphi},\quad {\rm div}({\bf \varphi} )={\rm div}_{x'}(\varphi')+\partial_{x_3}\varphi_3,\\
 \\
  \displaystyle  {\rm rot}({\bf \varphi} )= (\partial_{x_2}\varphi_3-\partial_{x_3}\varphi_2,
  -\partial_{x_1}\varphi_3+\partial_{x_3}\varphi_1,\partial_{x_1}\varphi_2-\partial_{x_2}\varphi_1)^t,
 \end{array}$$
where ${\bf \varphi} =(\varphi',\varphi_3)$ with $\varphi'=(\varphi_1, \varphi_2)$, is a vector  function defined in $\Omega^\epsilon$. \\

 Finally, we introduce some functional spaces. $L^{2}_0$  is the space of functions of $L^2$  with zero mean value. Let $C^\infty_{\#}(Z)$ be the space of infinitely differentiable functions in $\mathbb{R}^3$ that are $Z'$-periodic. By $L^2_{\#}(Z)$ (resp. $H^1_{\#}(Z)$),  we denote its completion in the norm $L^2(Z)$ (resp. $H^1(Z)$) and by $L^2_{0,\#}(Z)$  the space of functions in $L^2_{\#}(Z)$ with zero mean value.

\subsection{Micropolar equations with non-standard boundary conditions} In this paper, we consider the 3D stationary linearized micropolar equations in dimensionless form by setting 
$${\bf u}_\ep=({\bf u}'_\ep(x), u_{3,\ep}(x)),\quad {\bf w}_\ep=({\bf w}'_\ep(x), w_{3,\ep}(x)),\quad p_\ep=p_\ep(x),$$
at a point $x\in \Omega_\ep$, which is given by
\begin{equation}\label{system_1}
\left\{\begin{array}{rl}
\displaystyle -\Delta {\bf u}_\ep+\nabla p_\ep=2N^2{\rm rot}( {\bf w}_\ep) & \hbox{in}\ \Omega^\epsilon,\\
\\
{\rm div}( {\bf u}^\epsilon)=0& \hbox{in}\ \Omega^\epsilon,\\
\\
\displaystyle -R_M\Delta  {\bf w}_\ep+4N^2{\bf w}_\ep=2N^2{\rm rot}({\bf u}_\ep)& \hbox{in}\ \Omega^\epsilon.
\end{array}\right.
\end{equation}
Imposing that the mean value of the pressure is zero, we complete the above system with the following boundary conditions on the top 
\begin{equation}\label{BCBot}
\begin{array}{ll}
\displaystyle {\bf u}_\epsilon=0,\quad  {\bf w}_\epsilon=0&\hbox{on}\ \Gamma_1^\ep,
\end{array}
\end{equation}
the following boundary condition on the lateral boundary
\begin{equation}\label{BCBot_lat}
\displaystyle {\bf u}_\epsilon={\bf g}_\ep,\quad  {\bf w}_\epsilon=0\quad\hbox{on}\  \Gamma_{\rm lat}^\ep,\quad\hbox{with}\quad \int_{\Gamma_{\rm lat}^\ep}{\bf g}_\ep\cdot {\bf n}\,d\sigma=0,
\end{equation}
and the non-standard boundary conditions on the bottom  
\begin{equation}\label{BCTopBot1}
\begin{array}{rl}
\displaystyle {\bf w}_\ep\times {\bf n}={\alpha\over 2}{\rm rot}({\bf u}_\ep)\times {\bf n}&\hbox{on}\ \Gamma_0,\\
\noame
\displaystyle R_M{\rm rot}({\bf w}_\ep)\times {\bf n}=2N^2\beta({\bf u}_{\ep} - {\bf s})\times {\bf n}&\hbox{on}\ \Gamma_0,\\
\noame
\displaystyle {\bf u}_\ep\cdot {\bf n}=0,\quad {\bf w}_\ep\cdot {\bf n}=0&\hbox{on}\ \Gamma_0.
\end{array}
\end{equation}
Here, ${\bf s}=({\bf s}',0)$ with ${\bf s}'=(s_1,s_2)$, stands for the characteristic velocity of the moving surface. Dimensionless (non-Newtonian) parameter $N^2$ characterizes the coupling between the equations for the velocity and microrotation and it is of order $O(1)$. The second dimensionless parameter, denoted by $R_M$ is, in fact, related to the characteristic length of the microrotation effects and will be compared with small parameter $\ep$ (see condition \ref{RMnumber}). Parameter  $\alpha$ characterizes the microrotation on the solid surfaces (see \cite{Bayada_NewModel} for more details). The coefficient $\beta$ allows the control of the slippage at the wall when the value ${\bf u}_{\ep} -{\bf s}$ is not zero. Finally, the ${\bf n}$ is the outside normal vector to the boundary (observe that in the case of $\Gamma_0$ it holds ${\bf n}=-{\bf e}_3$).

In order to study the asymptotic behavior of the solution of system (\ref{system_1})--(\ref{BCTopBot1}), we also need to assume a certain regularity on the boundary data $g_\ep$, and uniform estimates of relevant norms. A very general way of stating those properties is the following: there exists a sequence of lift functions ${\bf J}_\ep\in H^1(\Omega_\ep)^2$ satisfying ${\rm div}({\bf J}_\ep)=0$ in $\Omega_\ep$, the boundary conditions
\begin{equation}\label{Jlift}
{\bf J}_\ep=0 \hbox{  on  }\Gamma_1^\ep,\quad {\rm J}_\ep={\bf g}_\ep \hbox{  on  }\Gamma_{\rm lat}^\ep,\quad {\bf J}_\ep\cdot {\bf n}=0 \hbox{  on  }\Gamma_0,
\end{equation}
and the estimates for every $\ep>0$
\begin{equation}\label{estimLift}
\|{\bf J}_\ep\|_{L^2(\Omega_\ep)^3}\leq C\ep^{1\over 2},\quad \|D{\bf J}_\ep\|_{L^2(\Omega_\ep)^{3\times 3}}\leq C\ep^{-{1\over 2}},\quad  \|{\bf J}_\ep\|_{L^2(\Gamma_0)^3}\leq C,
\end{equation}
where $C>0$ is a universal constant independent of $\ep$.

\begin{remark} One typical construction of a boundary data $g_\ep$ and the associate lift function ${\bf J}_\ep$ is the following, see \cite{Bayada_NewModel}. Consider a regular vector field ${\bf\widetilde J}\in H^1(\widetilde\Omega^\ep)^3$, satisfying
$${\rm div}({\bf \widetilde J})=0\ \ {\rm  in  }\ \ \widetilde\Omega_\ep,\quad  {\bf\widetilde J}=0\ \ {\rm  on  }\ \ \widetilde \Gamma_1^\ep,\quad {\bf\widetilde J}\cdot {\bf n}=0\ \ {\rm on  }\ \ \Gamma_0.$$
Extending ${\bf\widetilde J}$ by zero outside of $\widetilde\Omega^\ep$, we can define ${\bf J}_\ep\in H^1(\Omega^\ep)^3$ by 
$${\bf J}_\ep(x)=( {\bf \widetilde J}'(x',x_3/\ep)), \ep \widetilde J_3(x', x_3/\ep))\quad \forall \, x\in \widetilde \Omega_\ep.$$
Note that only  normal component of the velocity is known on $\Gamma_0$ and is zero, while the tangential velocity is not given. Nevertheless, we can  choose an artificial value ${\bf g}=({\bf g'},0)$, with ${\bf g}'=(g_1, g_2)$, of the velocity on $\Gamma_0$. We choose it in a way that function $\xi\in H^{1\over 2}(\partial\Omega^\ep)^3$ defined by
 $$\xi=\left\{\begin{array}{ll} 0&\hbox{on } \Gamma_1^\ep,\\
 \noame
 {\bf g}_\ep&\hbox{on } \Gamma_{\rm lat}^\ep,\\
 \noame
 ({\bf g}',0)&\hbox{on } \Gamma_0,
 \end{array}\right.$$
 satisfying  $\int_{\partial\Omega^\ep}\xi\cdot {\bf n}\,d\sigma=0$. \\

Then, by the change of variables $(x', x_3)=(x',\ep z_3)$, it holds
$$\begin{array}{l}
\displaystyle \int_{\Omega^\ep}|D{\bf   J}_\ep|^2\,dx=\ep\int_{\widetilde\Omega^\ep}\left(
|D_{x'}{\bf \widetilde J}'|^2+{1\over \ep^2}|\partial_{z_3}{\bf\widetilde J}'|^2+\ep^2|\nabla_{x'}\widetilde J_3|^2+|\partial_{z_3}\widetilde J_3|^2
\right)\,dx'dz_3,\\
\noame
\displaystyle
\int_{\Omega^\ep}|{\bf J}_\ep|^2\,dx=\ep\int_{\widetilde\Omega^\ep}(|{\bf \widetilde J}'|^2+\ep^2|\widetilde J_3|^2)\,dx'dz_3,
\end{array}$$
so that ${\bf J}_\ep$ satisfies all the required properties (\ref{Jlift}) and (\ref{estimLift}).

\end{remark}

Next, we are able to give the problem with homogeneous boundary conditions on the lateral and top boundaries. Introducing
$${\bf v}_\ep={\bf u}_\ep-{\bf J}_\ep,$$
system (\ref{system_1})--(\ref{BCTopBot1}) can be rewritten as follows
\begin{equation}\label{system_1_J}
\left\{\begin{array}{rl}
\displaystyle -\Delta ({\bf v}_\ep+{\bf J}_\ep)+\nabla p_\ep=2N^2{\rm rot}( {\bf w}_\ep) & \hbox{in}\ \Omega^\epsilon,\\
\\
{\rm div}( {\bf v}^\epsilon)=0& \hbox{in}\ \Omega^\epsilon,\\
\\
\displaystyle -R_M\Delta  {\bf w}_\ep+4N^2{\bf w}_\ep=2N^2\left({\rm rot}({\bf v}_\ep)+{\rm rot}({\bf J}_\ep)\right)& \hbox{in}\ \Omega^\epsilon,
\end{array}\right.
\end{equation}
 with the following boundary conditions 
\begin{equation}\label{BCTopBot1_J}
\begin{array}{rl}
\displaystyle {\bf v}_\epsilon=0,\quad  {\bf w}_\epsilon=0&\hbox{on}\ \Gamma_1^\ep\cup \Gamma_{\rm lat}^\ep,\\
\noame
\displaystyle {\bf w}_\ep\times {\bf n}={\alpha\over 2}\left({\rm rot}({\bf v}_\ep)+{\rm rot}({\bf J}_\ep)\right)\times {\bf n} &\hbox{on}\ \Gamma_0,\\
\noame
\displaystyle R_M{\rm rot}({\bf w}_\ep)\times {\bf n}=2N^2\beta({\bf v}_{\ep}+{\bf g} -{\bf s})\times {\bf n}&\hbox{on}\ \Gamma_0,\\
\noame
\displaystyle {\bf v}_\ep\cdot {\bf n}=0,\quad {\bf w}_\ep\cdot {\bf n}=0&\hbox{on}\ \Gamma_0.\\
\end{array}
\end{equation}
In the next subsection, we recall the existence and uniqueness result  of weak solution of problem (\ref{system_1_J})--(\ref{BCTopBot1_J}).
\subsection{Mathematical study}
In this paper, the goal is to derive an effective model describing  by using rigorous asymptotic analysis with respect to the small parameter $\ep$. In particular, we will focus on detecting the effects of the roughness together with the effects of non-standard boundary conditions. We remark that different asymptotic behaviors of the flow can be deduced depending on the order of magnitude of the dimensionless parameter. Indeed, if we compare the characteristic number $R_M$ with small parameter $\ep$, three different asymptotic situations can be formally identifed (see e.g. \cite{Bayada_Gamouana, BayadaLuc}). We consider the  most interesting one, which  leads to a strong
coupling at main order, namely the regime
\begin{equation}\label{RMnumber}R_M=\epsilon^2\,R_c,\quad R_c>0.
\end{equation}

Let us start by defining the notion of weak solution to system (\ref{system_1_J})--(\ref{BCTopBot1_J}). Due to the boundary conditions, we introduce the following functional spaces 
$$\begin{array}{c}
\displaystyle V^\ep=\{\varphi\in H^1(\Omega^\ep)^3\,:\, \varphi=0\quad \hbox{in }\Gamma_1^\ep\cup \Gamma_{\rm lat}^\ep,\quad \varphi\cdot {\bf n}=0\quad\hbox{ on }\Gamma_0\},\quad V_0^\ep=\{\varphi\in V^\ep\,:\, {\rm div}(\varphi)=0\},
\end{array}
$$
equipped with the norm of $\|D\varphi\|_{L^2(\Omega^\ep)^{3\times 3}}$,  and the space
$$L^2_0(\Omega^\ep)=\left\{q\in L^2(\Omega^\ep)\,:\, \int_{\Omega_\ep}q\,dx=0\right\},$$
equipped with the norm of $L^2(\Omega^\ep)$. Observe that it holds the following identities, which will be important in the following:
\begin{eqnarray}
&&\displaystyle -\Delta \varphi={\rm rot}({\rm rot} (\varphi))-\nabla {\rm div}(\varphi)\quad\forall \varphi\in \mathcal{D}(\Omega^\ep)^3,\label{Delta_Rot}\\
\noame
&&\displaystyle \int_{\Omega^\ep}{\rm rot}(\psi)\cdot \varphi\,dx=\int_{\Omega^\ep}{\rm rot}(\varphi)\cdot \psi\,dx-\int_{\Gamma_0}(\psi\times {\bf n})\cdot \varphi\,d\sigma(x')\quad\forall (\varphi,\psi)\in H^1(\Omega^\ep)^3\times H^1(\Omega^\ep)^3.\label{Parts}
\end{eqnarray}

\begin{proposition}[Theorem 2.1 in \cite{Bayada_NewModel}] Sufficiently regular solutions of  (\ref{system_1_J})--(\ref{BCTopBot1_J}) satisfy the weak formulation: 

Find $({\bf v}_\ep, w_\ep, p_\ep)\in V^\ep_0\times V^\ep\times L^2_0(\Omega^\ep)$ such that 
\begin{eqnarray}
&&\displaystyle \int_{\Omega_\ep}{\rm rot}({\bf v}_\ep)\cdot {\rm rot}(\varphi)\,dx-\int_{\Omega}p_\ep\,{\rm div}(\varphi)\,dx
-2N^2\int_{\Omega_\ep} {\rm rot}(\varphi)\cdot{\bf w}_\ep\,dx \label{Form_Var_vel}\\
\noame
&&\displaystyle +2\left(N^2-{1\over \alpha}\right)\int_{\Gamma_0}({\bf w}_\ep\times {\bf n})\cdot \varphi\,d\sigma(x')=-\int_{\Omega_\ep}{\rm rot}({\bf J}_\ep)\cdot {\rm rot}(\varphi)\,dx\quad\quad\forall \varphi\in V^\ep,\nonumber
\\
\noame
&&\displaystyle \ep^2 R_c\int_{\Omega_\ep}{\rm rot}({\bf w}_\ep)\cdot {\rm rot}(\psi)\,dx+  \ep^2 R_c\int_{\Omega_\ep}{\rm div}({\bf w}_\ep)\cdot {\rm div}(\psi)\,dx+4N^2\int_{\Omega_\ep}{\bf w}_\ep\cdot\psi\,dx\nonumber\\
\noame
&&\displaystyle  -2N^2\int_{\Omega_\ep}{\bf v}_\ep\cdot {\rm rot}(\psi)\,dx-2N^2(\beta-1)\int_{\Gamma_0}({\bf v}_\ep\times {\bf n})\cdot\psi\,d\sigma(x') =2N^2\int_{\Omega_\ep}{\bf J}_\ep \cdot {\rm rot}(\psi)\,dx\label{Form_Var_micro}\\
\noame
&&\displaystyle 
+2N^2\beta\int_{\Gamma_0}({\bf g}-{\bf s})\times {\bf n}\cdot \psi\,d\sigma(x')-2N^2\int_{\Gamma_0}({\bf g}\times{\bf n})\cdot\psi\,d\sigma(x')\quad\quad\forall \psi\in V^\ep.\nonumber
\end{eqnarray}
\end{proposition}
\begin{proof} The derivation of the variational formulations (\ref{Form_Var_vel}) and (\ref{Form_Var_micro}) is given in \cite[Theorem 2.1]{Bayada_NewModel}. Here, we reproduce some steps:
\begin{itemize}
\item[--]  Multiplying (\ref{system_1_J})$_1$ by $\varphi\in V^\ep$, using (\ref{Delta_Rot}) and (\ref{Parts}), we deduce
\begin{eqnarray}
&&\displaystyle \int_{\Omega_\ep}{\rm rot}({\bf v}_\ep)\cdot {\rm rot}(\varphi)\,dx-\int_{\Omega}p_\ep\,{\rm div}(\varphi)\,dx
-2N^2\int_{\Omega_\ep} {\rm rot}(\varphi)\cdot{\bf w}_\ep\,dx \label{Form_Var_vel_derivation1}\\
\noame
&&\displaystyle -\int_{\Gamma_0}({\rm rot}({\bf v}_\ep)\times {\bf n})\cdot \varphi\,d\sigma(x')-\int_{\Gamma_0}({\rm rot}({\bf J}_\ep)\times {\bf n})\cdot \varphi\,d\sigma(x')\nonumber\\
\noame
&&\displaystyle =-2N^2\int_{\Omega^\ep}({\bf w}_\ep\times {\bf n})\cdot \varphi\,d\sigma(x')
-\int_{\Omega^\ep}{\rm rot}({\bf J}_\ep)\cdot {\rm rot}(\varphi)\,dx.\nonumber
\end{eqnarray}
Taking into account the boundary condition (\ref{BCTopBot1_J})$_2$, we obtain (\ref{Form_Var_vel}).\\

\item[--] Multiplying (\ref{system_1_J})$_3$ by $\psi\in V^\ep$, using (\ref{Delta_Rot}) and (\ref{Parts}), we get
\begin{eqnarray}
&&\displaystyle  \ep^2 R_c\int_{\Omega_\ep}{\rm rot}({\bf w}_\ep)\cdot {\rm rot}(\psi)\,dx+\ep^2 R_c\int_{\Omega_\ep}{\rm div}({\bf w}_\ep)\cdot {\rm div}(\psi)\,dx+4N^2\int_{\Omega_\ep}{\bf w}_\ep\cdot\psi\,dx\nonumber\\
\noame
&&\displaystyle -\ep^2 R_c\int_{\Gamma_0}(({\rm rot}({\bf w}_\ep)\times {\bf n})\cdot \psi\,d\sigma(x')-2N^2\int_{\Omega_\ep}{\bf v}_\ep\cdot {\rm rot}(\psi)\,dx+2N^2\int_{\Gamma_0}({\bf v}_\ep\times {\bf n})\cdot\psi\,d\sigma(x')\label{Form_Var_micro_derivation1}
\\
\noame
 && =2N^2\int_{\Omega_\ep}{\bf J}_\ep \cdot {\rm rot}(\psi)\,dx-2N^2\int_{\Gamma_0}({\bf g}\times {\bf n})\cdot\psi\,d\sigma(x').\nonumber
\end{eqnarray}
Taking into account the boundary condition (\ref{BCTopBot1_J})$_3$, we obtain (\ref{Form_Var_micro}).
\end{itemize}
\end{proof}
To finish this section, we give the existence and uniqueness result, whos proof can be found in \cite[Theorem 2.2]{Bayada_NewModel} taking into account the rescaling of $R_M=\ep^2 R_c$ and that the maximum of the function $h_\ep$ is $\ep^2h_{\rm max}$.
\begin{theorem}
Under the previous assumptions and assuming that $\alpha$ and $\beta$ satisfy condition 
\begin{equation}\label{Existence_condition}
\gamma^2=\left|{1\over \alpha}-N^2-N^2\beta\right|^2<{R_c\over h_{\rm max}^2}(1-N^2),
\end{equation}
then, for every $\ep>0$,  there exists a unique weak solution $({\bf v}_\epsilon, {\bf w}_\epsilon, p_\epsilon)\in V^\ep_0\times V^\ep\times L^2_0(\Omega^\ep)$ of (\ref{Form_Var_vel})--(\ref{Form_Var_micro}).
\end{theorem}

\section{Rescaling}\label{sec:rescaling}
To study the asymptotic behavior of the solutions ${\bf u}_\ep$, ${\bf w}_\ep$ and $p_\ep$ when $\ep$  tends to zero,  we need to have the solutions in a domain with fixed height. To do this, we use the dilatation in the variable $x_3$ given by
\begin{equation}\label{dilatacion}
z_3={x_3\over \ep}\,,
\end{equation}
in order to have the functions defined in  $\widetilde\Omega_\ep$ and boundaries $\widetilde \Gamma_1^\epsilon$ and  $\widetilde \Gamma_{\rm lat}^\ep$ given in (\ref{domains_tilde}).\\

Using the change of variables (\ref{dilatacion}) in  (\ref{system_1_J})--(\ref{BCTopBot1_J}), we obtain the   rescaled system
\begin{equation}\label{system_1_dil}
\left\{\begin{array}{rl}
\displaystyle -\Delta_\ep( {\bf \widetilde v}_\ep+{\bf \widetilde J}_\ep)+\nabla_\ep \widetilde p_\ep=2N^2{\rm rot}_\ep( {\bf \widetilde w}_\ep) & \hbox{in}\ \widetilde \Omega_\epsilon,\\
\\
{\rm div}( {\bf \widetilde v}^\epsilon)=0& \hbox{in}\ \widetilde \Omega_\epsilon,\\
\\
\displaystyle -\ep^2 R_c\Delta_\ep   {\bf \widetilde w}_\ep+4N^2 {\bf \widetilde w}_\ep=2N^2\left({\rm rot}_\ep({\bf \widetilde u}_\ep)+{\rm rot}_\ep({\bf \widetilde J}_\ep)\right)& \hbox{in}\ \widetilde \Omega_\epsilon,
\end{array}\right.\end{equation}
with the rescaled boundary conditions
\begin{equation}\label{BCTopBot1_J_tilde}
\begin{array}{rl}
\displaystyle {\bf \widetilde v}_\epsilon=0,\quad  {\bf  \widetilde  w}_\epsilon=0&\hbox{on}\  \widetilde \Gamma_1^\ep\cup \widetilde \Gamma_{\rm lat}^\ep,\\
\noame
\displaystyle {\bf  \widetilde w}_\ep\times {\bf n}={\alpha\over 2}\left({\rm rot}_\ep({\bf  \widetilde v}_\ep)+{\rm rot}({\bf  \widetilde J}_\ep)\right)\times {\bf n} &\hbox{on}\ \Gamma_0,\\
\noame
\displaystyle \ep^2 R_c\,{\rm rot}_\ep({\bf \widetilde  w}_\ep)\times {\bf n}=2N^2\beta({\bf \widetilde  v}_{\ep}+{\bf g} -{\bf s})\times {\bf n}&\hbox{on}\ \Gamma_0,\\
\noame
\displaystyle {\bf \widetilde v}_\ep\cdot {\bf n}=0,\quad {\bf \widetilde w}_\ep\cdot {\bf n}=0&\hbox{on}\ \Gamma_0,
\end{array}
\end{equation}

\noindent where the unknown functions in the above system are given by ${\bf \widetilde u}_\epsilon(x',z_3)={\bf u}^\epsilon(x',\epsilon z_3)$, $ \widetilde p_\epsilon(x',z_3)=p_\epsilon(x',\epsilon z_3)$, ${\bf \widetilde w}_\epsilon(x',z_3)={\bf w}_\epsilon(x',\epsilon z_3)$ and ${\bf \widetilde J}_\epsilon(x',z_3)=(\widetilde {\bf J}'(x',z_3), \ep J_3(x',z_3))$  for a.e. $(x',z_3)\in \widetilde\Omega^\epsilon$. Here, for a vectorial function $\widetilde\varphi=(\widetilde\varphi',\widetilde\varphi_3)$ defined in $\widetilde\Omega_\ep$,  the rescaled operators are given by
\begin{equation}\label{def_operator_tilde}\begin{array}{c}
 \displaystyle \Delta_{\epsilon}{\bf \widetilde\varphi}=\Delta_{x'}{\bf \widetilde\varphi}+{1\over \ep^2}\partial_{z_3}^2{\bf \widetilde\varphi},\quad {\rm div}_{\epsilon}({\bf \widetilde \varphi})={\rm div}_{x'}(\widetilde \varphi')+{1\over \ep}\partial_{z_3}\widetilde \varphi_3,
 \\
 \noame
\displaystyle  {\rm rot}_{\epsilon}({\bf\widetilde  \varphi})=\left( {\rm rot}_{x'}(\widetilde\varphi_3)+{1\over \ep}{\rm rot}_{z_3}(\widetilde\varphi'), {\rm Rot}_{x'}(\widetilde\varphi')\right)^t,
 \end{array}
 \end{equation}
 where, denoting by 
 \begin{equation}\label{vec_perp}(\widetilde\varphi')^\perp=(-\widetilde\varphi_2,\widetilde\varphi_1)^t,
 \end{equation}
 we define
 \begin{equation}\label{rotationals}
 {\rm rot}_{x'}(\widetilde\varphi_3)=(\partial_{x_2}\widetilde\varphi_3,-\partial_{x_1}\widetilde\varphi_3)^t,\quad {\rm rot}_{z_3}(\widetilde\varphi')=\partial_{z_3}(\widetilde\varphi')^\perp,\quad {\rm Rot}_{x'}(\widetilde\varphi')=\partial_{x_1}\widetilde\varphi_2-\partial_{x_2}\widetilde\varphi_1.
 \end{equation}
 Moreover, from estimates (\ref{estimLift}), we also have
 \begin{equation}\label{Jlift_tilde}
 \|{\bf \widetilde J}_\ep\|_{L^2(\widetilde \Omega_\ep)^3}\leq C,\quad \|D_\ep{\bf\widetilde  J}_\ep\|_{L^2(\widetilde \Omega_\ep)^{3\times 3}}\leq C\ep^{-1},\quad \|{\bf\widetilde J}_\ep\|_{L^2(\Gamma_0)^3}\leq C.
 \end{equation}
In addition, the rescaled  variational formulation in $\widetilde \Omega_\ep$ is the following:

Find $({\bf \widetilde v}_\ep, {\bf \widetilde w}_\ep, \widetilde p_\ep)\in \widetilde V^\ep_0\times \widetilde V^\ep\times L^2_0(\widetilde \Omega_\ep)$ such that  
\begin{equation}\label{Form_Var_vel_tilde}
\begin{array}{l}
\displaystyle \ep\int_{\widetilde \Omega_\ep}{\rm rot}_\ep({\bf\widetilde  v}_\ep)\cdot {\rm rot}_\ep(\widetilde \varphi)\,dx'dz_3-\ep\int_{\widetilde \Omega_\ep}\widetilde p_\ep\,{\rm div}_\ep(\widetilde \varphi)\,dx'dz_3-2N^2\ep\int_{\widetilde \Omega_\ep} {\rm rot}_\ep(\widetilde \varphi)\cdot {\bf \widetilde w}_\ep\,dx'dz_3\\
\noame
\displaystyle +2\left(N^2-{1\over \alpha}\right)\int_{\Gamma_0}({\bf \widetilde w}_\ep\times {\bf n})\cdot\widetilde \varphi\,d\sigma(x') =-\ep\int_{\widetilde \Omega_\ep}{\rm rot}_\ep({\bf \widetilde J}_\ep)\cdot {\rm rot}_\ep(\widetilde \varphi)\,dx'dz_3,\\
\\
\displaystyle \ep^3R_c\int_{\widetilde \Omega_\ep}{\rm rot}_\ep( {\bf \widetilde w}_\ep)\cdot {\rm rot}_\ep(\widetilde \psi)\,dx'dz_3+\ep^3 R_c\int_{\widetilde \Omega_\ep}{\rm div}_\ep({\bf \widetilde w}_\ep)\cdot {\rm div}_\ep(\widetilde \psi)\,dx'dz_3+4N^2\ep\int_{\widetilde \Omega_\ep}{\bf \widetilde w}_\ep\cdot\widetilde \psi\,dx'dz_3\\
\noame
\displaystyle -2N^2\ep\int_{\widetilde\Omega_\ep}{\bf\widetilde v}_\ep\cdot {\rm rot}_\ep(\widetilde\psi)\,dx'dz_3

-2N^2(\beta-1)\int_{\Gamma_0}({\bf \widetilde v_\ep}\times {\bf n})\cdot\widetilde \psi\,d\sigma(x') \\
\noame
\displaystyle =2N^2\ep\int_{\widetilde \Omega_\ep} {\bf\widetilde  J}_\ep\cdot{\rm rot}_\ep(\widetilde \psi)\,dx'dz_3+2N^2\beta\int_{\Gamma_0}({\bf g}-{\bf s})\times {\bf n}\cdot \widetilde \psi\,d\sigma(x') 
-2N^2\int_{\Gamma_0}({\bf g}\times{\bf n})\cdot\widetilde \psi\,d\sigma(x'),
\end{array}
\end{equation}
for every $(\widetilde \varphi,\widetilde\psi)\in \widetilde V^\ep\times \widetilde V^\ep$ and $\widetilde q\in L^2(\widetilde \Omega_\ep)$ obtained from $(\varphi, \psi, q)$ by using the change of variables (\ref{dilatacion}), where 
$$\begin{array}{c}
\displaystyle \widetilde V^\ep=\{\widetilde \varphi\in H^1(\widetilde \Omega_\ep)^2\,:\, \widetilde \varphi=0\quad \hbox{in }\widetilde \Gamma^\ep_1\cup \widetilde \Gamma_{\rm lat}^\ep,\quad \widetilde \varphi\cdot {\bf n}=0\quad\hbox{ on }\Gamma_0\},\quad \widetilde V_0^\ep=\{\widetilde \varphi\in \widetilde V^\ep\,:\, {\rm div}_\ep(\widetilde \varphi)=0\}.
\end{array}
$$
In the next section, we establish the {\it a priori estimates} for the rescaled unknown in an $\ep$-independent domain.

\section{{\it A priori} estimates and convergences} \label{sec:estimates}
 This section is devoted to derive the {\it a priori} estimates of the unknowns and is divided in three parts. First, we deduce  the {\it a priori} estimates for velocity and microrotation in an $\ep$-independent domain, i.e. the domain $\Omega$ defined in Section \ref{sec:setting}. Next, we extend the pressure to $\Omega$ and derive its corresponding {\it a priori} estimates. Finally, from previos {\it a priori} estimates, we deduce some convergence results. 

\subsection{Estimates for velocity and microroration} To derive the desired estimates,  let us recall and proof some well-known technical results (see for instance \cite{Anguiano_SG}).
\begin{lemma}[Poincar\'e's  inequality]\label{Poincare_lemma} For all $\varphi\in V^\ep$, there holds the following inequality
\begin{equation}\label{Poincare}
\|\varphi\|_{L^2(\Omega^\epsilon)^3}\leq C\epsilon\|D \varphi\|_{L^2(\Omega^\epsilon)^{3\times 3}}.
\end{equation}
Moreover, from the change of variables (\ref{dilatacion}),   there holds for all $\widetilde\varphi\in \widetilde V^\ep$ the following rescaled estimate
\begin{equation}\label{Poincare2}
\|\widetilde \varphi\|_{L^2(\widetilde \Omega^\epsilon)^3}\leq C\epsilon\|D_{\ep} \widetilde \varphi\|_{L^2(\widetilde \Omega^\epsilon)^{3\times 3}}.
\end{equation}
\end{lemma}
\begin{lemma}
For all $\varphi \in V^\ep$, it holds
\begin{equation}\label{Gaffney}
 \|D\varphi\|_{L^2(\Omega^\ep)^{3\times 3}}^2=\|{\rm rot}(\varphi)\|^2_{L^2(\Omega^\ep)^{3}}+\|{\rm div}(\varphi)\|^2_{L^2(\Omega^\ep)}, \quad \|{\rm rot}(\varphi)\|^2_{L^2(\Omega^\ep)^{3}}\leq \|D\varphi\|_{L^2(\Omega^\ep)^{3\times 3}}^2,
\end{equation}
and, if moreover, ${\rm div}(\varphi)=0$ in $\Omega^\ep$, then 
\begin{equation}\label{Gaffney_div0}
\|{\rm rot}(\varphi)\|_{L^2(\Omega_\ep)^{3}}^2=\|D\varphi\|_{L^2(\Omega_\ep)^{3\times 3}}^2.
\end{equation}
Moreover, from the change of variables (\ref{dilatacion}), we have for all $\widetilde\varphi\in \widetilde V^\ep$ that 
\begin{equation}\label{Gaffney_tilde}\|D_\ep \widetilde \varphi\|^2_{L^2(\widetilde \Omega^\ep)^{3\times 3}}=\|{\rm rot}_\ep(\widetilde \varphi)\|^2_{L^2(\widetilde \Omega^\ep)^{3}}+\|{\rm div}_\ep(\widetilde \varphi)\|^2_{L^2(\widetilde \Omega^\ep)},\quad \|{\rm rot}_\ep(\widetilde \varphi)\|^2_{L^2(\widetilde \Omega^\ep)^{3}}\leq \|D_\ep \widetilde \varphi\|^2_{L^2(\widetilde \Omega^\ep)^{3\times 3}},
\end{equation}
and if moreover, ${\rm div}_\ep(\widetilde \varphi)=0$ in $\widetilde \Omega^\ep$, then  
\begin{equation}\label{Gaffney_div0_tilde}
\|{\rm rot}_\ep\widetilde \varphi\|_{L^2(\widetilde \Omega_\ep)^{3}}^2=\|D_\ep\widetilde \varphi\|_{L^2(\widetilde \Omega_\ep)^{3\times 3}}^2.
\end{equation}
\end{lemma}
\begin{proof} For all $\varphi\in V^\ep$, it holds (see for instance \cite{Boyer} formula (IV.23))
$$\int_{\Omega^\ep}(|{\rm div}(\varphi)|^2+|{\rm rot}(\varphi)|^2)\,dx=\int_{\Omega^\ep}|D\varphi|^2\,dx+\int_{\Gamma_0}((\varphi\cdot\nabla){\bf n})\cdot \varphi\,d\sigma(x').$$
Since $\Gamma_0$ is flat, then the term $\int_{\Gamma_0}((\varphi\cdot\nabla){\bf n})\cdot \varphi\,d\sigma(x')$ vanishes, ans so we get (\ref{Gaffney}). Then, (\ref{Gaffney_div0}) is a consequence of the free divergence condition.

Finally, the change of variables (\ref{dilatacion}) applied to (\ref{Gaffney}) and (\ref{Gaffney_div0}) implies (\ref{Gaffney_tilde}) and (\ref{Gaffney_div0_tilde}), respectively, for any $\widetilde\varphi$ obtained from $\varphi$ by the change of variables (\ref{dilatacion}).

\end{proof}

\begin{lemma}[Trace estimates] For all $\widetilde\varphi \in \widetilde V^\ep$, the following inequalities hold\begin{equation}\label{trace_estimate}
 \|\widetilde\varphi'\|_{L^2(\Gamma_0)^2}\leq \ep h_{\rm max}^{1\over 2}\|D_\ep \widetilde\varphi\|_{L^2(\widetilde \Omega^\ep)^{3\times 3}},\quad \|\widetilde\varphi\times {\bf n}\|_{L^2(\Gamma_0)^2}\leq \ep h_{\rm max}^{1\over 2}\|D_\ep \widetilde\varphi\|_{L^2(\widetilde \Omega^\ep)^{3\times 3}}.
\end{equation}
\end{lemma}
\begin{proof}
Thank to $\widetilde \varphi(x',h(x'/\ep^\ell))=0$ in $\omega$, we have that
$$\begin{array}{rl}
\displaystyle \int_{\Gamma_0}|\widetilde \varphi|^2d\sigma= &\displaystyle \int_{\omega}|\widetilde \varphi(x',0)|^2\,dx'= \int_{\omega}\left|\int_0^{h(x'/\ep^\ell)}\partial_{z_3}\widetilde \varphi(x',z_3)\,dz_3\right|^2\,dx'\leq   h_{\rm max}\int_{\widetilde \Omega^\ep}|\partial_{z_3}\widetilde \varphi|^2\,dx'dz_3,
\end{array}$$
that is, 
$$\|\widetilde \varphi'\|_{L^2(\Gamma_0)^2}\leq h_{\rm max}^{1\over 2}\|\partial_{z_3}\widetilde \varphi'\|_{L^2(\widetilde \Omega^\ep)^{2}},$$
which implies (\ref{trace_estimate})$_1$. Since
$$\|\widetilde\varphi\times {\bf n}\|_{L^2(\Gamma_0)^3}=\|(-\widetilde\varphi_2, \widetilde\varphi_1, 0)\|_{L^2(\Gamma_0)^3},$$ 
from (\ref{trace_estimate})$_1$, we deduce (\ref{trace_estimate})$_2$.

\end{proof}

\begin{lemma}[A priori estimates]\label{lemma_estimates} Under hypothesis (\ref{Existence_condition}), the following estimates hold for the solution $({\bf \widetilde v}_\ep, {\bf\widetilde w}_\ep)$ of problem (\ref{Form_Var_vel_tilde}):
\begin{eqnarray}
\displaystyle
\|{\bf \widetilde v}^\epsilon\|_{L^2(\widetilde \Omega^\epsilon)^3}\leq C, &\displaystyle
\|D_{\epsilon} {\bf  \widetilde v}^\epsilon\|_{L^2(\widetilde\Omega^\epsilon)^{3\times 3}}\leq C\epsilon^{-1},& \label{estim_sol_dil1}\medskip
\\
\displaystyle
\|{\bf \widetilde w}^\epsilon\|_{L^2(\widetilde\Omega^\epsilon)^3}\leq C\ep^{-1}, &\displaystyle
\|D_{\epsilon} {\bf \widetilde w}^\epsilon\|_{L^2(\widetilde\Omega^\epsilon)^{3\times 3}}\leq C\epsilon^{-2}.\label{estim_sol_dil2}
\end{eqnarray}
As consequence, it also holds
\begin{eqnarray}
\displaystyle
\|{\bf \widetilde u}^\epsilon\|_{L^2(\widetilde \Omega^\epsilon)^3}\leq C, &\displaystyle
\|D_{\epsilon} {\bf \widetilde u}^\epsilon\|_{L^2(\widetilde\Omega^\epsilon)^{3\times 3}}\leq C\epsilon^{-1}.& \label{estim_sol_dil1_u}
\end{eqnarray}
\end{lemma}
\begin{proof}  The proof of (\ref{estim_sol_dil1}) and (\ref{estim_sol_dil2})  is similar to the one given in \cite[Lemma 3.1]{Bayada_NewModel}, but here we are in the 3D case, so we will give some remarks of the development. We remark that, once proved estimates  (\ref{estim_sol_dil1}),  the proof of (\ref{estim_sol_dil1_u})  is straightforward from
$${\bf \widetilde u}_\ep={\bf \widetilde v}_\ep+{\bf \widetilde J}_\ep,$$
and  estimates of ${\bf \widetilde J}_\ep$ given by (\ref{Jlift_tilde}). We will prove (\ref{estim_sol_dil1}) and (\ref{estim_sol_dil2}) in 4 steps. \\

{\it Step 1. Mixed variational formulation.} Problem (\ref{Form_Var_vel_tilde}) can be written in the following mixed variational form:

Find $({\bf  \widetilde v}_\ep, {\bf  \widetilde w}_\ep,  \widetilde p_\ep)\in  \widetilde V^\ep\times \widetilde V^\ep\times L^2_0(\widetilde \Omega^\ep)$, such that
\begin{eqnarray}
\mathcal{A}^\ep(({\bf  \widetilde v}_\ep, {\bf  \widetilde w}_\ep), ( \widetilde \varphi,\widetilde \psi))+\mathcal{B}^\ep((\widetilde \varphi,\widetilde \psi),\widetilde p_\ep)=\mathcal{L}(\widetilde \varphi,\widetilde \psi)&\forall\,(\widetilde \varphi,\widetilde \psi)\in  \widetilde V^\ep\times \widetilde V^\ep,\label{Form_Var_A}\\
\noame
\mathcal{B}^\ep(({\bf  \widetilde v}_\ep, {\bf  \widetilde w}_\ep),\widetilde q)=0&\forall\, \widetilde q\in L^2_0(\widetilde\Omega^\ep),
\end{eqnarray}
with $\mathcal{A}^\ep$ defined by
\begin{equation}\label{A_bilinear_1}
\begin{array}{rl}
\mathcal{A}^\ep(({\bf  \widetilde v}_\ep, {\bf  \widetilde w}_\ep), ( \widetilde \varphi,\widetilde \psi))=& \displaystyle \ep\int_{\widetilde \Omega_\ep}{\rm rot}_\ep({\bf\widetilde  v}_\ep)\cdot {\rm rot}_\ep(\widetilde \varphi)\,dx'dz_3-2N^2\ep\int_{\widetilde \Omega_\ep} {\rm rot}_\ep(\widetilde \varphi)\cdot {\bf \widetilde w}_\ep\,dx'dz_3\\
\noame
&\displaystyle +\ep^3R_c\int_{\widetilde \Omega_\ep}{\rm rot}_\ep( {\bf \widetilde w}_\ep)\cdot {\rm rot}_\ep(\widetilde \psi)\,dx'dz_3+\ep^3 R_c\int_{\widetilde \Omega_\ep}{\rm div}_\ep({\bf \widetilde w}_\ep)\cdot {\rm div}_\ep(\widetilde \psi)\,dx'dz_3\\
\noame
&\displaystyle +4N^2\ep\int_{\widetilde \Omega_\ep}{\bf \widetilde w}_\ep\cdot\widetilde \psi\,dx'dz_3
  -2N^2\ep\int_{\widetilde\Omega_\ep}{\bf\widetilde v}_\ep\cdot {\rm rot}_\ep(\widetilde\psi)\,dx'dz_3
\\
\noame
&\displaystyle +2\left(N^2-{1\over \alpha}\right)\int_{\Gamma_0}({\bf \widetilde w}_\ep\times {\bf n})\cdot\widetilde \varphi\,d\sigma(x')-2N^2(\beta-1)\int_{\Gamma_0}({\bf \widetilde v_\ep}\times {\bf n})\cdot\widetilde \psi\,d\sigma(x'),
\end{array}
\end{equation}
with $\mathcal{B}^\ep$ given by
\begin{equation}\label{B_linear}
\mathcal{B}^\ep((\widetilde \varphi,\widetilde \psi),\widetilde p_\ep)=-\ep\int_{\widetilde\Omega^\ep}\widetilde p_\ep\,{\rm div}_\ep(\widetilde\varphi)\,dx'dz_3,
\end{equation}
and $\mathcal{L}^\ep$ given by
\begin{equation}\label{L_linear}
\displaystyle\mathcal{L}^\ep(\widetilde \varphi,\widetilde \psi)= 2N^2\ep\int_{\widetilde \Omega_\ep} {\rm rot}_\ep({\bf\widetilde  J}_\ep)\cdot \widetilde \psi\,dx'dz_3-\ep\int_{\widetilde \Omega_\ep}{\rm rot}_\ep({\bf \widetilde J}_\ep)\cdot {\rm rot}_\ep(\widetilde \varphi)\,dx'dz_3+2N^2\beta\int_{\Gamma_0}({\bf g}-{\bf s})\times {\bf n}\cdot \widetilde \psi\,d\sigma(x').
\end{equation}
Taking into account that 
$$2\left(N^2-{1\over \alpha}\right)\int_{\Gamma_0}({\bf \widetilde w}_\ep\times {\bf n})\cdot\widetilde \varphi\,d\sigma=-2\left(N^2-{1\over \alpha}\right)\int_{\Gamma_0}(\widetilde \varphi\times {\bf n})\cdot{\bf \widetilde w}_\ep\,d\sigma(x'),$$
and 
$$\begin{array}{rl}
\displaystyle-2N^2\ep\int_{\widetilde\Omega_\ep}{\bf\widetilde v}_\ep\cdot {\rm rot}_\ep(\widetilde\psi)\,dx'dz_3=&\displaystyle 
-2N^2\ep\int_{\widetilde\Omega_\ep}\widetilde\psi\cdot {\rm rot}_\ep({\bf\widetilde v}_\ep)\,dx'dz_3+2N^2\int_{\Gamma_0}(\widetilde\psi\times {\bf n})\cdot {\bf\widetilde v}_\ep\,d\sigma(x')\\
\noame
=&
\displaystyle 
-2N^2\ep\int_{\widetilde\Omega_\ep}\widetilde\psi\cdot {\rm rot}_\ep({\bf\widetilde v}_\ep)\,dx'dz_3-2N^2\int_{\Gamma_0}({\bf\widetilde v}_\ep\times {\bf n})\cdot \widetilde\psi\,d\sigma(x'),
\end{array}$$
we then give another expression for $\mathcal{A}^\ep$ by
\begin{equation}\label{A_bilinear}
\begin{array}{rl}
\mathcal{A}^\ep(({\bf  \widetilde v}_\ep, {\bf  \widetilde w}_\ep), ( \widetilde \varphi,\widetilde \psi))=&  \displaystyle \ep\int_{\widetilde \Omega_\ep}{\rm rot}_\ep({\bf\widetilde  v}_\ep)\cdot {\rm rot}_\ep(\widetilde \varphi)\,dx'dz_3-2N^2\ep\int_{\widetilde \Omega_\ep} {\rm rot}_\ep(\widetilde \varphi)\cdot {\bf \widetilde w}_\ep\,dx'dz_3\\
\noame
&\displaystyle +\ep^3R_c\int_{\widetilde \Omega_\ep}{\rm rot}_\ep( {\bf \widetilde w}_\ep)\cdot {\rm rot}_\ep(\widetilde \psi)\,dx'dz_3+\ep^3 R_c\int_{\widetilde \Omega_\ep}{\rm div}_\ep({\bf \widetilde w}_\ep)\cdot {\rm div}_\ep(\widetilde \psi)\,dx'dz_3\\
\noame
&\displaystyle +4N^2\ep\int_{\widetilde \Omega_\ep}{\bf \widetilde w}_\ep\cdot\widetilde \psi\,dx'dz_3
 -2N^2\ep\int_{\widetilde\Omega_\ep}\widetilde\psi\cdot {\rm rot}_\ep({\bf\widetilde v}_\ep)\,dx'dz_3
\\
\noame
&\displaystyle -2\left(N^2-{1\over \alpha}\right)\int_{\Gamma_0}(\widetilde \varphi\times {\bf n})\cdot{\bf \widetilde w}_\ep\,d\sigma(x')-2N^2\beta\int_{\Gamma_0}({\bf \widetilde v_\ep}\times {\bf n})\cdot\widetilde \psi\,d\sigma(x').
\end{array}
\end{equation}

{\it Step 2. Estimate $\mathcal{A}^\ep(({\bf \widetilde v}_\ep, {\bf \widetilde w}_\ep), ({\bf \widetilde v}_\ep, {\bf \widetilde w}_\ep))$ from below.}
From identities (\ref{Gaffney_tilde}) applied to ${\bf \widetilde v}_\ep$ and (\ref{Gaffney_div0_tilde}) applied to ${\bf \widetilde w}_\ep$, applying Cauchy-Schwarz's inequality and using estimates (\ref{trace_estimate}), we deduce
$$\begin{array}{rl}
\mathcal{A}^\ep(({\bf \widetilde v}_\ep, {\bf \widetilde w}_\ep), ({\bf \widetilde v}_\ep, {\bf \widetilde w}_\ep))=&\displaystyle \ep\int_{\widetilde \Omega_\ep}|{\rm rot}_\ep({\bf\widetilde  v}_\ep)|^2\,dx'dz_3-4N^2\ep\int_{\widetilde \Omega_\ep} {\rm rot}_\ep({\bf\widetilde  v}_\ep)\cdot {\bf \widetilde w}_\ep\,dx'dz_3\\
\noame
&\displaystyle +\ep^3R_c\int_{\widetilde \Omega_\ep}|{\rm rot}_\ep( {\bf \widetilde w}_\ep)|^2\,dx'dz_3+\ep^3 R_c\int_{\widetilde \Omega_\ep}|{\rm div}_\ep({\bf \widetilde w}_\ep)|^2\,dx'dz_3\\
\noame
&\displaystyle +4N^2\ep\int_{\widetilde \Omega_\ep}|{\bf \widetilde w}_\ep|^2 \,dx'dz_3 +2\gamma \int_{\Gamma_0}({\bf\widetilde  v}_\ep\times {\bf n})\cdot{\bf \widetilde w}_\ep\,d\sigma\\
\noame
\geq &\ep\|D_\ep{\bf \widetilde v}_\ep\|^2_{L^2(\widetilde\Omega^\ep)^{3\times 3}}+\ep^3R_c\|D_\ep{\bf \widetilde w}_\ep\|^2_{L^2(\widetilde\Omega^\ep)^{3\times 3}}+4N^2\ep\|{\bf \widetilde w}_\ep\|^2_{L^2(\widetilde\Omega^\ep)^{3}}
\\
\noame
&\displaystyle -4N^2\ep\|D_\ep{\bf \widetilde v}_\ep\|_{L^2(\widetilde\Omega^\ep)^{3\times 3}}\| {\bf \widetilde w}_\ep\|_{L^2(\widetilde\Omega^\ep)^{3}}-2|\gamma|\ep^2 h_{\rm max}\|D_\ep{\bf \widetilde v}_\ep\|_{L^2(\widetilde\Omega^\ep)^{3\times 3}}\|D_\ep{\bf \widetilde w}_\ep\|_{L^2(\widetilde\Omega^\ep)^{3\times 3}}.\end{array}$$
 From condition (\ref{Existence_condition}), there exists $c_1>0$ such that
\begin{equation}\label{c1cond}
{\gamma h_{\rm max}\over R_c}<c_1<{1-N^2\over \gamma h_{\rm max}}.
\end{equation}
Applying Young's inequality
\begin{equation}\label{c1condYoung}
\|D_\ep{\bf\widetilde v_\ep}\|_{L^2(\widetilde\Omega_\ep)^{3\times 3}}\|D_\ep{\bf\widetilde w_\ep}\|_{L^2(\widetilde\Omega_\ep)^{3\times 3}}\leq {c_1\over 2\ep}\|D_\ep{\bf\widetilde v_\ep}\|_{L^2(\widetilde\Omega_\ep)^{3\times 3}}^2+{\ep\over 2c_1}\|D_\ep{\bf\widetilde w_\ep}\|^2 _{L^2(\widetilde\Omega_\ep)^{3\times 3}}.
\end{equation}
By continuity, there exists a real number $c_2$ satisfying $0<c_2<c_1$, and such that
\begin{equation}\label{c2cond}
c_1<{1-{N^2\over c_2}\over \gamma h_{\rm max}},
\end{equation}
and applying again Young's inequality, we also have
\begin{equation}\label{cccondYoung}
\|D_\ep{\bf\widetilde v_\ep}\|_{L^2(\widetilde\Omega_\ep)^{3\times 3}}\| {\bf\widetilde w_\ep}\|_{L^2(\widetilde\Omega_\ep)^{3 }}\leq {1\over 4c_2}\|D_\ep{\bf\widetilde v_\ep}\|^2_{L^2(\widetilde\Omega_\ep)^{3\times 3}}+{c_2}\| {\bf\widetilde w_\ep}\|^2_{L^2(\widetilde\Omega_\ep)^{3}}.
\end{equation}
Then, we deduce
\begin{equation}\label{Acoerciveness}\begin{array}{rl}
\mathcal{A}^\ep(({\bf \widetilde v}_\ep, {\bf \widetilde w}_\ep), ({\bf \widetilde v}_\ep, {\bf \widetilde w}_\ep))\geq &
\displaystyle
\ep\left(1-{N^2\over c_2}-|\gamma| h_{\rm max} c_1\right)\|D_\ep{\bf\widetilde  v}_\ep\|^2_{L^2(\widetilde\Omega^\ep)^{3\times 3}}\\
\noame
 &\displaystyle+\ep^3\left(
  R_c- {|\gamma|h_{\rm max}\over c_1}\right)\|D_\ep{\bf\widetilde  w}_\ep\|^2_{L^2(\widetilde\Omega^\ep)^{3\times 3}}+4N^2(1-c_2)\ep\| {\bf\widetilde  w}_\ep\|^2_{L^2(\widetilde\Omega^\ep)^{3 }}\\
 \noame
 \geq&\displaystyle  \ep A\|D_\ep{\bf\widetilde  v}_\ep\|^2_{L^2(\widetilde\Omega^\ep)^{3\times 3}}+\ep^3B\|D_\ep{\bf\widetilde  w}_\ep\|^2_{L^2(\widetilde\Omega^\ep)^{3\times 3}},
\end{array}
\end{equation}
with
$$A=\left(1-{N^2\over c_2}-|\gamma| h_{\rm max} c_1\right),\quad B=\left(
 R_c- {|\gamma|h_{\rm max}\over c_1}\right).$$
 Using conditions (\ref{c1cond}) and (\ref{c2cond}), $A$ and $B$ are positive.\\

 {\it Step 3. Estimate $\mathcal{L}^\ep({\bf \widetilde v}_\ep, {\bf \widetilde w}_\ep)$.} From Cauchy-Schwarz's inequality (written here with constant $C_P$) and (\ref{Gaffney}), we get
 \begin{equation}\label{L_linear_estim}
\begin{array}{rl}
\displaystyle|\mathcal{L}^\ep({\bf \widetilde v}_\ep, {\bf \widetilde w}_\ep)|=&\displaystyle \left|2N^2\ep\int_{\widetilde \Omega^\ep} {\rm rot}_\ep({\bf\widetilde  J}_\ep)\cdot {\bf \widetilde w}_\ep\,dx'dz_3 -\ep\int_{\widetilde \Omega_\ep}{\rm rot}_\ep({\bf \widetilde J}_\ep)\cdot {\rm rot}_\ep({\bf \widetilde v}_\ep)\,dx'dz_3+2N^2\beta\int_{\Gamma_0}({\bf g}-{\bf s})\times {\bf n}\cdot {\bf \widetilde w}_\ep\,d\sigma\right|\\
\noame
\leq &\displaystyle 2N^2\ep^2C_P\|{\rm rot}_\ep{\bf \widetilde J}_\ep\|_{L^2(\widetilde\Omega^\ep)^{3\times 3}}\|D_\ep{\bf \widetilde w}_\ep\|_{L^2(\widetilde\Omega^\ep)^{3\times 3}}+\ep\|D_\ep{\bf \widetilde J}_\ep\|_{L^2(\widetilde\Omega^\ep)^{3\times 3}}\|D_\ep{\bf \widetilde v}_\ep\|_{L^2(\widetilde\Omega^\ep)^{3\times 3}}\\
\noame
&\displaystyle +2N^2\beta\ep h_{\rm max}^{1\over 2} \|{\bf g}-{\bf s}\|_{L^2(\Gamma_0)^2}\|D_\ep{\bf \widetilde w}_\ep\|_{L^2(\widetilde\Omega^\ep)^{3\times 3}}\\
\noame
=&\displaystyle \ep\|D_\ep{\bf \widetilde J}_\ep\|_{L^2(\widetilde\Omega^\ep)^{3\times 3}}\|D_\ep{\bf \widetilde v}_\ep\|_{L^2(\widetilde\Omega^\ep)^{3\times 3}}+2N^2\ep C_\ep\|D_\ep{\bf \widetilde w}_\ep\|_{L^2(\widetilde\Omega^\ep)^{3\times 3}},
\end{array}
\end{equation}
where $C_\ep=\ep C_P\|{\rm rot}_\ep{\bf \widetilde J}_\ep\|_{L^2(\widetilde\Omega^\ep)^{3\times 3}}
+\beta h^{1\over 2}_{\rm max} \|{\bf g}-{\bf s}\|_{L^2(\Gamma_0)^2}$. Observe that from properties (\ref{Jlift_tilde}) satisfied by ${\bf J}_\ep$, there exists a a constant $C>0$ such that $C_\ep\leq C$ for any $\ep>0$.\\

 {\it Step 4. Derivation of the a priori estimates.} To obtain estimates (\ref{estim_sol_dil1}) and (\ref{estim_sol_dil2}), we use inequalities (\ref{Acoerciveness}) and (\ref{L_linear_estim}) to get
\begin{equation}\label{estim_together}
\begin{array}{l}
\displaystyle  \ep A\|D_\ep{\bf\widetilde  v}_\ep\|^2_{L^2(\widetilde\Omega^\ep)^{3\times 3}}+\ep^3B\|D_\ep{\bf\widetilde  w}_\ep\|^2_{L^2(\widetilde\Omega^\ep)^{3\times 3}}\\
\noame
\leq \displaystyle \ep\|D_\ep{\bf \widetilde J}_\ep\|_{L^2(\widetilde\Omega^\ep)^{3\times 3}}\|D_\ep{\bf \widetilde v}_\ep\|_{L^2(\widetilde\Omega^\ep)^{3\times 3}}+2N^2\ep C_\ep\|D_\ep{\bf \widetilde w}_\ep\|_{L^2(\widetilde\Omega^\ep)^{3\times 3}},
\end{array}
\end{equation}
which is equivalent to the following inequality
\begin{equation}\label{estim_together2}
\begin{array}{l}
\displaystyle  \ep A\|D_\ep{\bf\widetilde  v}_\ep\|^2_{L^2(\widetilde\Omega^\ep)^{3\times 3}}-\ep\|D_\ep{\bf \widetilde J}_\ep\|_{L^2(\widetilde\Omega^\ep)^{3\times 3}}\|D_\ep{\bf \widetilde v}_\ep\|_{L^2(\widetilde\Omega^\ep)^{3\times 3}}\\
\noame
\displaystyle +\ep^3B\|D_\ep{\bf\widetilde  w}_\ep\|^2_{L^2(\widetilde\Omega^\ep)^{3\times 3}}-2N^2\ep C_\ep\|D_\ep{\bf \widetilde w}_\ep\|_{L^2(\widetilde\Omega^\ep)^{3\times 3}}\leq 0.
\end{array}
\end{equation}
Summing in both sides ${\ep\over 4 A}\|D_\ep{\bf\widetilde  J}_\ep\|^2_{L^2(\widetilde\Omega^\ep)^{3\times 3}}$ and ${N^4\over \ep B}C^2_\ep$, we observe that last inequality can be written equivalently as
\begin{equation}\label{estim_together3}
\begin{array}{l}
\displaystyle  \left(\ep^{1\over 2} A^{1\over 2}\|D_\ep{\bf\widetilde  v}_\ep\|_{L^2(\widetilde\Omega^\ep)^{3\times 3}}-{\ep^{1\over 2}\over 2A^{1\over 2}}\|D_\ep{\bf\widetilde  J}_\ep\|_{L^2(\widetilde\Omega^\ep)^{3\times 3}}\right)^2+\left(\ep^{3\over 2}B^{1\over 2}\|D_\ep{\bf\widetilde  w}_\ep\|_{L^2(\widetilde\Omega^\ep)^{3\times 3}}-{N^2\over \ep^{1\over 2} B^{1\over 2}} C_\ep\right)^2\\
\noame
\leq \displaystyle {\ep\over 4 A}\|D_\ep{\bf\widetilde  J}_\ep\|^2_{L^2(\widetilde\Omega^\ep)^{3\times 3}}+{N^4\over \ep B}C^2_\ep.
\end{array}
\end{equation}
which is written as follows
\begin{equation}\label{estim_together4}
\begin{array}{l}
\displaystyle  \ep A\left(\|D_\ep{\bf\widetilde  v}_\ep\|_{L^2(\widetilde\Omega^\ep)^{3\times 3}}-{1\over 2A}\|D_\ep{\bf\widetilde  J}_\ep\|_{L^2(\widetilde\Omega^\ep)^{3\times 3}}\right)^2+\ep^3B\left(\|D_\ep{\bf\widetilde  w}_\ep\|_{L^2(\widetilde\Omega^\ep)^{3\times 3}}-{N^2\over \ep^2 B} C_\ep\right)^2\\
\noame
\leq \displaystyle {\ep\over 4 A}\|D_\ep{\bf\widetilde  J}_\ep\|^2_{L^2(\widetilde\Omega^\ep)^{3\times 3}}+{N^4\over \ep B}C^2_\ep.
\end{array}
\end{equation}
Since $C_\ep$ is  uniformly bounded, using assumptions (\ref{Jlift_tilde}), we deduce than the right hand side of the previous inequality is bounded by $C/\ep$  for a certain constant $C>0$. This implies (\ref{estim_sol_dil1})$_2$ and (\ref{estim_sol_dil2})$_2$. Finally, from the Poincar\'e inequality (\ref{Poincare2}), we get estimates (\ref{estim_sol_dil1})$_1$ and (\ref{estim_sol_dil2})$_1$, respectively. This finishes the proof.

 \end{proof}

\paragraph{The extension of $({\bf \widetilde u}^\epsilon, {\bf \widetilde w}^\epsilon)$  to the whole domain.}
The sequence of solutions $({\bf \widetilde u}^\ep, {\bf \widetilde w}^\ep)$ is defined in a varying set $\widetilde\Omega^\ep$, but not defined in the fixed domain $\Omega$ independent of $\ep$. Thus, to pass to the limit if $\ep$ tends to zero, we need convergences in fixed Sobolev spaces (defined in $\Omega$), so we need to extend  $({\bf \widetilde  u}^\ep, {\bf \widetilde w}^\ep)$ to the whole domain $\Omega$. From the boundary conditions satisfied by ${\bf \widetilde u}^\epsilon$ and ${\bf\widetilde w}^{\epsilon}$ on the top boundary $\widetilde \Gamma_1^\ep$,  we can extend them  by zero in $\Omega\setminus \widetilde{\Omega}^{\epsilon}$ and we denote the extensions by   the same
symbol.
\begin{lemma}[Estimates of extended unknowns]\label{lemma_estimates2} The extended functions $({\bf \widetilde u}^\ep, {\bf \widetilde w}^\ep)$ satisfy the following estimates
\begin{eqnarray}
\displaystyle
\|{\bf \widetilde u}_\epsilon\|_{L^2(  \Omega )^3}\leq C ,& \displaystyle
\|D_{\epsilon} {\bf \widetilde u}_\epsilon\|_{L^2( \Omega )^{3\times 3}}\leq C\epsilon^{-1}, \label{estim_sol_dil_ext_u}
\medskip
\\
\displaystyle
\|{\bf \widetilde w}_\epsilon\|_{L^2( \Omega)^3}\leq C\ep^{-1}, & \displaystyle
\|D_{\epsilon}{\bf  \widetilde w}_\epsilon\|_{L^2( \Omega )^{3\times 3}}\leq C\epsilon^{-2}. \label{estim_sol_dil_ext_w}
\end{eqnarray}

\end{lemma}
\begin{proof}  Estimates for the extension of ${\bf \widetilde u}_\epsilon$ and ${\bf \widetilde w}_\epsilon$ are obtained straightforward from  (\ref{estim_sol_dil1}) and (\ref{estim_sol_dil2}), respectively. 

\end{proof}

\subsection{Estimates for pressure}

Extending the pressure is a much more difficult task.  A continuation of the pressure for a flow in a porous media was introduced in \cite{Tartar}. This construction applies to periodic holes in a domain $\Omega_\epsilon$ when each hole is strictly contained into the periodic cell. In this context, we can not use directly this result because the ``holes" are along the oscillating boundary $\Gamma_1^\ep$ of $\Omega_\epsilon$ and moreover, the scale of the vertical direction is smaller than the scales of the horizontal directions. This fact will induce several limitations in the results obtained by using the method, especially in view of the convergence for the pressure. In this sense, for the case of Newtonian fluids,  an operator $R^\epsilon$  generalizing the results of \cite{Tartar}  to this context (extending the pressure from $\Omega_\ep$ to $Q_\ep$) was introduced in \cite{Bayada_Chambat, Mikelic2}, and later extended to the case of non-Newtonian (power law) fluids in \cite{Anguiano_SG} by defining an extension operator $R^\ep_p$, for every flow index $p>1$. 

Then,  to extend the pressure to the whole domain $\Omega$, the mapping $R^\ep$  (defined in  \cite[Lemma 4.6]{Anguiano_SG} as $R_2^\ep$) allows us to extend the pressure $p_\ep$ from $\Omega_\ep$ to $Q_\ep$ by introducing $F_\ep$ in $H^{-1}(Q_\ep)^3$ as follows (brackets are for duality products between $H^{-1}$ and $H^1_0$)
\begin{equation}\label{F}\langle F_\epsilon, \varphi\rangle_{Q_\epsilon}=\langle \nabla p_\epsilon, R^\epsilon (\varphi)\rangle_{\Omega_\epsilon}\quad \hbox{for any }\varphi\in H^1_0(Q_\epsilon)^3\,.
\end{equation}
We compute the right-hand side of (\ref{F}) by using the first equation of (\ref{Form_Var_vel}), and taking into account ${\bf u}_\ep={\bf v}_\ep+{\bf J}_\ep$,  we have
\begin{equation}\label{equality_duality}
\begin{array}{rl}
\displaystyle
\left\langle F_{\epsilon},\varphi\right\rangle_{Q_\epsilon}=&\displaystyle
 -\int_{\Omega_\ep}{\rm rot}({\bf u}_\ep)\cdot {\rm rot}(R^{\epsilon}(\varphi))\,dx
 + 2N^2 \int_{\Omega_\ep}{\rm rot}\,w_\ep\cdot R^\ep(\varphi)\,dx\,.
\end{array}\end{equation}
Using Lemma \ref{lemma_estimates} for fixed $\epsilon$, we see that it is a bounded functional on $H^1_0(Q_\epsilon)$ (see the proof of Lemma \ref{lemma_est_P} below) and in fact $F_\epsilon\in H^{-1}(Q_\epsilon)^3$. Moreover, ${\rm div} \varphi=0$ implies $\left\langle F_{\epsilon},\varphi\right\rangle_{Q_\epsilon}=0$, 
and the DeRham theorem gives the existence of $P_\epsilon$ in $L^{2}_0(Q_\epsilon)$ with $F_\epsilon=\nabla P_\epsilon$.

Defining the rescaled extended pressure $\widetilde P_\ep\in L^2_0(\Omega)$ by 
$$\widetilde P_\ep(x',z_3)=P_\ep(x',\ep z_3),\quad\hbox{a.e. }(x',z_3)\in\Omega,$$
we get for any $\widetilde \varphi\in H^1_0(\Omega)^3$ where $\widetilde\varphi(x',z_3)=\varphi(x',\ep z_3)$ that 
$$\begin{array}{rl}\displaystyle\langle \nabla_{\epsilon}\widetilde P_\epsilon, \widetilde\varphi\rangle_{\Omega}&\displaystyle
=-\int_{\Omega}\widetilde P_\epsilon\,{\rm div}_{\epsilon}\,\widetilde\varphi\,dx'dz_3
=- \epsilon^{-1}\int_{Q_\epsilon}P_\epsilon\,{\rm div}\,\varphi\,dx= \epsilon^{-1}\langle \nabla P_\epsilon, \varphi\rangle_{Q_\epsilon}\,.
\end{array}$$
Then, using the identification (\ref{equality_duality}) of $F_\epsilon$, we get
$$\begin{array}{rl}\displaystyle\langle \nabla_{\epsilon}\widetilde P_\epsilon, \widetilde\varphi\rangle_{\Omega}&\displaystyle
=\epsilon^{-1}\Big(-\int_{\Omega_\epsilon}{\rm rot}({\bf u}_\ep) : {\rm rot}(R^{\epsilon}(\varphi))\,dx
+2N^2\int_{\Omega_\ep}{\rm rot}({\bf w}_\ep)\cdot R^\ep(\varphi)\,dx\Big)\,,
\end{array}$$
and applying the change of variables (\ref{dilatacion}), we obtain 
\begin{equation}\label{extension_1}
\begin{array}{rl}\displaystyle\langle \nabla_{\epsilon}\widetilde P_\epsilon, \widetilde\varphi\rangle_{\Omega}
=&\displaystyle- \int_{\widetilde \Omega_\epsilon}{\rm rot}_{\ep}( {\bf \widetilde u}_\ep) : {\rm rot}_{\ep}( \widetilde R^{\epsilon}(\widetilde \varphi))\,dx'dz_3+2N^2\int_{\widetilde \Omega^\ep}{\rm rot}_{\ep}({\bf \widetilde w}_\ep)\cdot \widetilde R^\ep(\widetilde \varphi)\,dx'dz_3\,,
\end{array}
\end{equation}
where $\widetilde R^\ep(\widetilde\varphi)=R^\ep(\varphi)$ for any $ \varphi\in H^1_0(Q_\epsilon)^3$ where $\widetilde\varphi(x',z_3)=\varphi(x', \ep z_3)$.\\

Now, we estimate the right-hand side of (\ref{extension_1}) to obtain the {\it a priori} estimate of the pressure $\widetilde P_\ep$.

\begin{lemma}[Estimates for extended pressure]\label{lemma_est_P}
Under hypothesis (\ref{Existence_condition}), the following estimates hold for the the extension $\widetilde P_\ep\in L^2_0(\Omega)$ of the pressure $\widetilde p_\ep$
\begin{equation}\label{esti_P}
\|\widetilde P_\ep\|_{L^2(\Omega)}\leq C\ep^{-2},\quad \|\nabla_{\ep}\widetilde P_\ep\|_{H^{-1}(\Omega)^3}\leq C\ep^{-2}.
\end{equation}
\end{lemma}
\begin{proof} From the proof of Lemma 4.7-(i) in \cite{Anguiano_SG}, we have that $\widetilde R^\ep(\widetilde\varphi)$ satisfies the following estimates
\begin{equation}\label{ext_1}
\begin{array}{l}
\displaystyle\|\widetilde R^\ep(\widetilde\varphi)\|_{L^2(\widetilde\Omega_\ep)^3}\leq C\left(\|\widetilde\varphi\|_{L^2(\Omega)^3}+
\ep^\ell\sum_{j=1}^2\|\partial_{x_i}\widetilde\varphi\|_{L^2(\Omega)^{3}}  + \|\partial_{z_3}\widetilde\varphi\|_{L^2(\Omega)^3}\right)\,,\\
\noame
\displaystyle
\|\partial_{x_i}\widetilde R^\ep(\widetilde\varphi)\|_{L^2(\widetilde\Omega_\ep)^{3\times 2}}\leq  C\left({1\over \ep^\ell}\|\widetilde\varphi\|_{L^2(\Omega)^3}+
\sum_{j=1}^2\|\partial_{x_i}\widetilde\varphi\|_{L^2(\Omega)^{3}} + {1\over \ep^\ell}\|\partial_{z_3}\widetilde\varphi\|_{L^2(\Omega)^3}\right)\,,\quad i=1, 2,\\
\noame
\displaystyle\|\partial_{z_3}\widetilde R^\ep(\widetilde\varphi)\|_{L^2(\widetilde\Omega_\ep)^3}\leq C\left(\|\widetilde\varphi\|_{L^2(\Omega)^3}+
\ep^\ell\sum_{j=1}^2\|\partial_{x_i}\widetilde\varphi\|_{L^2(\Omega)^{3}}  + \|\partial_{z_3}\widetilde\varphi\|_{L^2(\Omega)^3}\right)\,.
\end{array}
\end{equation}
Thus, since $0<\ep< \ep^\ell\ll 1$, we deduce 
\begin{equation}\label{ext_2}
\|\widetilde R^\ep(\widetilde \varphi)\|_{L^2(\widetilde\Omega_\ep)^3}\leq C\|\widetilde \varphi\|_{H^1_0(\Omega)^3},\quad \|{\rm rot}_{\ep}(\widetilde R^\ep(\widetilde \varphi))\|_{L^2(\widetilde\Omega_\ep)^{3\times 3}}\leq {C\over \ep}\|\widetilde \varphi\|_{H^1_0(\Omega)^3}\,.
\end{equation}
By using Cauchy-Schwarz's inequality, estimate (\ref{Gaffney}), estimates for $D_{\ep}{\bf \widetilde u}_\ep$ in (\ref{estim_sol_dil1_u}) and for $D_{\ep}{\bf \widetilde w}_\ep$ in (\ref{estim_sol_dil2}) together with (\ref{ext_2}), we  obtain
\begin{equation}\label{estim_nabla_p_crit_1}
\begin{array}{rl}
\displaystyle
\left|\int_{\widetilde\Omega_\ep}{\rm rot}_{\ep}({\bf \widetilde u}_\ep):{\rm rot}_{\ep}(\widetilde R^\ep(\widetilde\varphi))\,dx'dz_3\right|\leq &\displaystyle  C\ep^{-1}\|D_{\ep}\widetilde R^\ep(\widetilde\varphi)\|_{L^2(\widetilde \Omega_\ep)^{3\times 3}}\leq  C\ep^{-2}\|\widetilde\varphi\|_{H^1_0(\Omega)^3},\\
\noame
\displaystyle \left|\int_{\widetilde\Omega_\ep}{\rm rot}_{\ep}({\bf \widetilde w}_\ep)\cdot \widetilde R^\ep(\widetilde\varphi)\,dx'dz_3\right|\leq &
\displaystyle \|D_{\ep}{\bf \widetilde w}_\ep\|_{L^2(\widetilde\Omega_\ep)^{3\times 3}}\|\widetilde R^\ep(\widetilde\varphi)\|_{L^2(\widetilde\Omega_\ep)^3} \\
\noame
\leq &  \displaystyle C\ep^{-2}\|\widetilde R^\ep(\widetilde\varphi)\|_{L^2(\widetilde\Omega_\ep)^3}C\|\widetilde\varphi\|_{H^1_0(\Omega)^3},
\end{array}
\end{equation}
which together with (\ref{extension_1}) gives $\|\nabla_{\ep}\widetilde P_\ep\|_{H^{-1}(\Omega)^3}\leq C\ep^{-2}$, i.e. estimate (\ref{esti_P})$_2$. By using the Ne${\breve{\rm c}}$as inequality there exists a representative $\widetilde P_\ep\in L^2_0(\Omega)$ such that
$$
\|\widetilde P_\ep\|_{L^2(\Omega)}\leq C\|\nabla\widetilde P_\ep\|_{H^{-1}(\Omega)^3}\leq C\|\nabla_{\ep}\widetilde P_\ep\|_{H^{-1}(\Omega)^3},
$$
which implies (\ref{esti_P})$_1$.

\end{proof}

\subsection{Convergences}
We give a compactness result concerning the asymptotic behavior of the extended sequences $({\bf \widetilde u}_\ep, {\bf \widetilde w}_\ep, \widetilde P_\ep)$  satisfying the {\it a priori} estimates given in Lemmas \ref{lemma_estimates2} and \ref{lemma_est_P}.
\begin{lemma}\label{lem_asymp_sub}
For a subsequence of $\ep$ still denote by $\ep$, we have the following convergence results:
\begin{itemize}
\item[(i)] {\it (Velocity)} There exists ${\bf \widetilde u}=({\bf \widetilde u'},\widetilde u_{3})\in H^1(0,h_{\rm max};L^2(\omega)^3)$, with $\widetilde u_{3}\equiv 0$ and ${\bf \widetilde u}'=0$ on $ \Gamma_1$, such that
\begin{eqnarray}
&\displaystyle {\bf  \widetilde u}^\ep\rightharpoonup {\bf \widetilde u}\hbox{  in  }H^1(0,h_{\rm max};L^2(\omega)^3),\label{conv_u_sub_tilde}\\
\noame
&\displaystyle {\rm div}_{x'}\left(\int_0^{h_{\rm max}}{\bf \widetilde u}'(x',z_3)\,dz_3\right)=0\quad \hbox{  in  }\omega, &\label{div_x_sub_tilde}\\
\noame
&\displaystyle \left(\int_0^{h_{\rm max}}{\bf \widetilde u}'(x',z_3)\,dz_3\right)\cdot {\bf n}'=0\quad \hbox{  on  }\partial\omega.&\nonumber
 \end{eqnarray}
\item[(ii)] {\it (Microrotation)} There exist ${\bf \widetilde w}=({\bf \widetilde w}',\widetilde w_3)\in H^1(0,h_{\rm max};L^2(\omega)^3)$, with $\widetilde w_3\equiv 0$ and ${\bf \widetilde w}'=0$ on $ \Gamma_1$, such that
\begin{eqnarray}
&\displaystyle \ep {\bf \widetilde w}_\ep\rightharpoonup {\bf\widetilde w}\hbox{  in  }H^1(0,h_{\rm max};L^2(\omega)^3).\label{conv_w_sub_tilde}
\end{eqnarray}
\item[(iii)] {\it (Pressure)} There  exist a function  $\widetilde p\in L^2_0(\omega)\cap H^1(\omega)$ (i.e. $\widetilde p$ is independent of $z_3$ and has mean value zero in $\omega$),  such that
\begin{eqnarray}
&\displaystyle \epsilon^{2}  \widetilde P_\ep\to  \widetilde p\hbox{  in  }L^2(\Omega).& \label{conv_P_sub}
\end{eqnarray}
\end{itemize}
\end{lemma}
\begin{proof}We start proving $(i)$. We will only give some remarks and,  for more details,  we refer the reader to Lemmas 5.2-i) and 5.4-i) in \cite{Anguiano_SG}.

We start with the extension $\tilde u_\ep$. Estimates  (\ref{estim_sol_dil_ext_u}) imply the existence of $\tilde u\in H^1(0,h_{\rm max};L^2(\omega)^3)$ such that convergence (\ref{conv_u_sub_tilde}) holds, and the continuity of the trace applications from the space of $\tilde u$ such that $\|{\bf \widetilde u}\|_{L^2}$ and $\|\partial_{z_3}{\bf \widetilde u}\|_{L^2}$ are bounded to $L^2(\Gamma_1)$ implies $\tilde u=0$ on $\Gamma_1$. It can also be deduced that $\widetilde u_3=0$ on $\Gamma_0$ according to boundary condition ${\bf \widetilde u}_\ep\cdot {\bf n}=\widetilde u_{3,\ep}=0$ on $\Gamma_0$.

Next, from the free divergence condition ${\rm div}_{\ep}({\bf \widetilde u}_\ep)=0$, it can be deduced that $\widetilde u_3$ is independent of $z_3$, which together with the boundary conditions satisfied by $\widetilde u_3$ on $\Gamma_1$ implies that $\widetilde u_3\equiv 0$. Finally, from the free divergence condition and the convergence (\ref{conv_u_sub_tilde}) of ${\bf \widetilde u}_\ep$, it is straightforward the corresponding free divergence condition in a thin domain given in (\ref{div_x_sub_tilde}).\\

We continue proving $(ii)$.  From estimates (\ref{estim_sol_dil_ext_w}), the  convergence of (\ref{conv_w_sub_tilde}) and that ${\bf \widetilde w}=0$ on $\Gamma_1$ and $\widetilde w_3=0$ on $\Gamma_0$ are obtained straighfordward.   It remains to prove that $\tilde w_3\equiv0$. To do this, we consider as test function $\widetilde \psi(x',z_3)=(0,0,\widetilde \psi_3)$ with $\widetilde\psi_3\in \mathcal{D}(\Omega)$,  in the variational formulation (\ref{Form_Var_vel_tilde}) extended to $\Omega$, taking into account ${\bf \widetilde u}_\ep={\bf \widetilde v}_\ep+{\bf \widetilde J}_\ep$ and the definition of the rotational given in (\ref{rotationals}), we get
$$
\begin{array}{l}
\displaystyle \ep^3 R_c\int_{\Omega}{\rm rot}_{x'}(\widetilde w_{3,\ep})\cdot {\rm rot}_{x'}(\widetilde \psi_3)\,dx'dz_3+\ep R_c\int_{\Omega}\partial_{z_3}\widetilde w_{3,\ep}\,\partial_{z_3}\widetilde \psi_3\,dx'dz_3+4N^2\ep\int_{\Omega}{\widetilde w}_{3,\ep}\, \widetilde \psi_3\,dx'dz_3\\
\noame
\displaystyle -2N^2\ep\int_{\Omega}{\bf\widetilde u}'_\ep\cdot {\rm rot}_{x'}(\widetilde\psi_3)\,dx'dz_3 =0.
\end{array}
$$
Passing to the limit by using convergences of $\widetilde u_\ep$ and $\widetilde w_\ep$ given respectively in (\ref{conv_u_sub_tilde}) and (\ref{conv_w_sub_tilde}), we get
$$R_c\int_\Omega \partial_{z_3}\widetilde w_3\,\partial_{z_3}\widetilde \psi_3\,dx'dz_3+4N^2\int_\Omega \widetilde w_3\,\widetilde \psi_3\,dx'dz_3=0\,,$$
and taking into account that $\widetilde w_3=0$ on $\Gamma_1\cup \Gamma_0$,  it is easily deduced that $\widetilde w_3\equiv 0$  in $\Omega$.\\

We finish the proof with $(iii)$. Estimate (\ref{esti_P}) implies, up to a subsequence, the existence of $\tilde P\in L^2_0(\Omega)$ such that 
\begin{equation}\label{conv_P_sub}
\ep^2 \widetilde P_\ep\rightharpoonup \tilde p\quad \hbox{  in  }L^2(\Omega).
\end{equation} Also, from $\|\nabla_{\ep}\widetilde P_\ep\|_{L^2(\Omega)^3}\leq C\ep^2$, by noting that $\partial_{z_3}\widetilde P_\ep/\ep^3$ also converges weakly in $H^{-1}(\Omega)$, we obtain $\partial_{z_3}\widetilde P=0$ and so $\widetilde P$ is independent of $z_3$. Since $\widetilde P_\ep$ has mean value zero in $\Omega$, then $\widetilde p$ also has mean value zero. To prove strong convergence, we refer to \cite[Lemma 7]{Anguiano_SG}. Moreover, since $\widetilde P_\ep$ has mean value zero in $\Omega$ and $\widetilde p$ does not depend on $z_3$, it holds
$$
  \int_\Omega \ep^{2}\widetilde P_\ep\,dx'dz_3\to \int_\Omega\widetilde p\,dx'dz_3=h_{\rm max}\int_{\omega}\widetilde p(x')\,dx'=0,
$$
which implies $\widetilde p\in L^2_0(\omega)$. Finally, the proof of $\widetilde p\in H^1(\omega)$ can be found in \cite[Theorem 3.1]{Bayada_NewModel}, so we omit it.

\end{proof}

\section{Adaptation of the unfolding method}\label{sec:unfolding}
\subsection{Definition and properties} 
As seen in Lemma \ref{lem_asymp_sub}, the change of variables (\ref{dilatacion}) does not capture the oscillations of the domain $\widetilde\Omega^\ep$. To capture them,  we use an adaptation of the unfolding method (see \cite{Ciora, Ciora2} for more details) introduced to this context in \cite{Anguiano_SG}.\\

\begin{remark}  For $k'\in \mathbb{Z}^2$, we define $\kappa: \mathbb{R}^2\to \mathbb{Z}^2$ by
\begin{equation}\label{kappa_fun}
\kappa(x')=k' \Longleftrightarrow x'\in Z_{k',1}\,.
\end{equation}
Remark that $\kappa$ is well defined up to a set of zero measure in $\mathbb{R}^2$ (the set $\cup_{k'\in \mathbb{Z}^2}\partial Z_{k',1}$). Moreover, for every $\epsilon>0$, we have
$$\kappa\left({x'\over \epsilon^\ell}\right)=k'\Longleftrightarrow x'\in Z_{k',\epsilon}\,.$$
\end{remark}

 \noindent Let us recall that this adaptation of the unfolding method divides the domain $\widetilde\Omega_\ep$ in cubes of lateral length $\ep^\ell$ and vertical length $h(z')$. Thus, given  $\widetilde\varphi\in L^2(\widetilde\Omega^\ep)^3$ (assuming $\widetilde\varphi$ extended by zero outside of $\omega$), we define $\widehat \varphi_\ep\in L^2(\mathbb{R}^2\times Z)^3$ by
\begin{equation}\label{varphi_hat}
{\bf \widehat\varphi}_{\epsilon}(x',z)={\bf \widetilde\varphi}\left( {\epsilon}^\ell\kappa\left(\frac{x' }{{\epsilon}^\ell} \right)+{\epsilon}^\ell z' ,z_3 \right)\text{\ \ a.e. \ }(x',z)\in \omega\times Z,
\end{equation}
where  the function $\kappa$ is defined by (\ref{kappa_fun}). Also, given $\widetilde\psi\in L^2(\Omega)$, we define $\widehat \psi_\ep\in L^2(\mathbb{R}^2\times \Pi)$ by\\
\begin{equation}\label{psi_hat}
{\bf \widehat\psi}_{\epsilon}(x',z)={\bf \widetilde\psi}\left( {\epsilon}^\ell\kappa\left(\frac{x' }{{\epsilon}^\ell} \right)+{\epsilon}^\ell z' ,z_3 \right)\text{\ \ a.e. \ }(x',z)\in \omega\times \Pi.
\end{equation}
  \begin{remark}\label{rem:hat} For $k'\in T_\ep$, the restriction of $\widehat \varphi_\ep$ to $Z'_{k',\ep}\times Z$ and $\widehat\psi_\ep$ to $Z'_{k',\ep}\times \Pi$ do not depend on $x'$, while as a function of $z$ it is  obtained from $(\widetilde \varphi,  \widetilde \psi)$  by using the change of variables
  $$z'={x'-\ep^\ell k'\over \ep^\ell},$$
  which transform $Z_{k',\ep}$ into $Z$ and $\widetilde Q_{k',\ep}$ into $\Pi$, respectively.
  \end{remark}
Below, we give some properties of the change of variables (\ref{varphi_hat}), whose can be found in  \cite[Lemma 4.9]{Anguiano_SG} in the case $p=2$.
\begin{proposition}\label{properties_hat} We have the following properties concerning the estimates of vectorial functions $\widetilde \varphi$  and scalar functions $\psi$ and their respective unfolding function $\widehat \varphi_\ep$  and $\widehat \psi_\ep$ given by (\ref{varphi_hat}) and (\ref{psi_hat}):
\begin{itemize}
\item[i)] For every $\widetilde\varphi\in L^2(\widetilde\Omega^\ep)^3$, we have
\begin{equation}\label{hat_varphi_prop1}
\|\widehat \varphi_\ep\|_{L^2(\omega\times Z)^3}=\|\widetilde \varphi\|_{L^2(\widetilde\Omega^\ep)^3}.
\end{equation}
\item[ii)] For every $\widetilde\varphi\in H^1(\widetilde\Omega^\ep)^3$, we have
\begin{equation}\label{hat_varphi_prop2}
\|D_{z'}\widehat \varphi_\ep\|_{L^2(\omega\times Z)^{3\times 3}}=\ep^\ell\|D_{x'} \widetilde \varphi\|_{L^2(\widetilde\Omega^\ep)^{3\times 3}},\quad \|\partial_{z_3}\widehat \varphi_\ep\|_{L^2(\omega\times Z)^{3}}=\|\partial_{z_3}\widetilde \varphi\|_{L^2(\widetilde\Omega^\ep)^{3}}.
\end{equation}
\item[iii)] For every $\widetilde\psi\in L^2(\Omega)$,   we have
\begin{equation}\label{hat_varphi_prop1}
\|\widehat \psi_\ep\|_{L^2(\omega\times \Pi)}=\|\widetilde \psi\|_{L^2(\Omega)}.
\end{equation}
\end{itemize}
\end{proposition}

In a similar way, let us introduce the adaption of the unfolding method on the boundary  $\Gamma_0$  (see Cioranescu et al. \cite{Ciora2} for more details). For this purpose, given $\widetilde \varphi\in L^2(\Gamma_0)^3$, we define $\widehat \varphi^b_\ep\in L^2(\mathbb{R}^2\times \widehat \Gamma_0)^3$ by
\begin{equation}\label{varphi_hat_b}
{\bf \widehat \varphi}_{\epsilon}^b(x',z)={\bf \widetilde \varphi} \left( {\epsilon}^\ell\kappa\left(\frac{x' }{{\epsilon}^\ell} \right)+{\epsilon}^\ell z' ,z_3 \right)\text{\ \ a.e. \ }(x',z)\in \omega\times \widehat \Gamma_0,
\end{equation}
where  the function $\kappa$ is defined by (\ref{kappa_fun}).

\begin{remark}\begin{itemize}
\item[(i)] Observe that from this definition, if we consider $\widetilde\varphi\in L^2(\Gamma_0)$, a $Z'$-periodic function, we define $\widetilde \varphi_\ep(x',z_3)=\widetilde\varphi(x'/\ep^\ell,z_3)$, it follows $\widehat \varphi^b_\ep(x',z)=\widetilde\varphi(z)$.

\item[(ii)] Observe that for $\widetilde\varphi\in H^1(\Gamma_0)$, $\widehat \varphi^b_\ep$ is the trace on $\Gamma_0$ of $\widehat \varphi_\ep$. Therefore, $\widehat \varphi^b_\ep$ has similar properties as $\widehat \varphi_\ep$. So it holds the following property:
\begin{equation}\label{hat_varphi_prop_b}
\|\widehat \varphi_\ep^b\|_{L^2(\omega\times \widehat \Gamma_0)^3}= \|\widetilde \varphi_\ep\|_{L^2(\Gamma_0)^3}.
\end{equation}
\end{itemize}
\end{remark}

\begin{definition}We define the unfolded unknowns as follows:
\begin{itemize}
\item[--] The unfolded velocity ${\bf \widehat{u}}_{\epsilon}$ is defined by applying (\ref{varphi_hat}) for $\widetilde\varphi={\bf \widetilde u}_\ep$.
\item[--] The unfolded microrotation ${\bf \widehat {w}}_{\epsilon}$ is defined by applying (\ref{varphi_hat}) for $\widetilde\varphi={\bf \widetilde w}_\ep$.
\item[--] The unfolded pressure $\widehat P_{\epsilon}$ is defined by applying (\ref{psi_hat}) for $\widetilde\psi=\widetilde P_\ep$.
\end{itemize}
\end{definition}

\noindent We are now in position to obtain estimates for the unfolded unknowns $({\bf \widehat{u}}_{\epsilon}, {\bf \widehat w}_\ep, \widehat P_\epsilon)$.
 \begin{lemma}\label{estimates_hat}
 There exists a constant $C>0$ independent of $\ep$, such that ${\bf \widehat u}_\ep$, ${\bf \widehat w}_\ep$ and $\widehat P_\epsilon$   satisfy
 \begin{eqnarray}
 \|{\bf \widehat u}_\ep\|_{L^2(\omega\times Z)^3}\leq C,&
 \|D_{z'}{\bf \widehat u}_\ep\|_{L^2(\omega\times Z)^{3\times 2}}\leq C\epsilon^{\ell-1},&
 \|\partial_{z_3}{\bf \widehat u}_\ep\|_{L^2(\omega\times Z)^{3}}\leq C,\label{estim_u_hat}\\
 \noame
  \|{\bf \widehat w}_\ep\|_{L^2(\omega\times Z)^3}\leq C\ep^{-1} ,&
 \|D_{z'}{\bf \widehat w}_\ep\|_{L^2(\omega\times Z)^{3\times 2}}\leq C\epsilon^{\ell-2},&
 \|\partial_{z_3}{\bf \widehat w}_\ep\|_{L^2(\omega\times Z)^3}\leq C\ep^{-1},\label{estim_w_hat}\\
 \noame
  &\|\widehat P_\ep\|_{L^2(\omega\times \Pi)}\leq C\ep^{-2}.&\label{estim_P_hat}
\label{estim_P01_hat}
  \end{eqnarray}
 \end{lemma}
\begin{proof}
The result is a consequence of the combination of  estimates given in Proposition \ref{properties_hat} with estimates for ${\bf \widetilde u}_\ep$, ${\bf \widetilde w}_\ep$,  $\widetilde P_\epsilon$ given in Lemma \ref{lemma_estimates} and  \ref{lemma_est_P}.

\end{proof}

\subsection{Equivalent weak variational formulation} We give the equivalent weak variational formulation of system  (\ref{Form_Var_vel_tilde}), which will be useful in next sections in order to obtain the limit system taking into account the effects of the rough boundary.

\noindent We consider $\widetilde \varphi_\ep(x',z_3)=\ep \widehat \varphi(x',x'/\ep^\ell,z_3)$ and $\widetilde\psi_\ep=\widehat \psi(x',x'/\ep^\ell,z_3)$ as test function in (\ref{Form_Var_vel_tilde}), where $\widehat \varphi, \widehat \psi\in \mathcal{D}(\omega;C^\infty_\#(Z)^3)$, and taking into account the extension of the pressure, we have 
$$
\langle\nabla_\ep \widetilde p_\ep,\widetilde\varphi_\ep\rangle_{\widetilde\Omega^\ep}=\langle\nabla_\ep \widetilde P_\ep\cdot \widetilde\varphi_\ep\rangle_{\Omega},
$$
and taking into account ${\bf \widetilde u}_\ep={\bf \widetilde v}_\ep+{\bf \widetilde J}_\ep$, we get
\begin{eqnarray}
&&\displaystyle \ep\int_{\widetilde \Omega_\ep}{\rm rot}_\ep({\bf\widetilde  u}_\ep)\cdot {\rm rot}_\ep(\widetilde\varphi_\ep)\,dx'dz_3-\ep\int_{\Omega}\widetilde P_\ep\,{\rm div}_\ep(\widetilde\varphi_\ep)\,dx'dz_3-2N^2\ep\int_{\widetilde \Omega_\ep} {\rm rot}_\ep(\widetilde\varphi_\ep)\cdot {\bf \widetilde w}_\ep\,dx'dz_3\label{Form_Var_vel_tilde_hat}\\
\noame
&&\displaystyle +2\left(N^2-{1\over \alpha}\right)\int_{\Gamma_0}({\bf \widetilde w}_\ep\times {\bf n})\cdot\widetilde \varphi_\ep\,d\sigma(x') =0,\nonumber
\\
\noame
&&\displaystyle \ep^3R_c\int_{\widetilde \Omega_\ep}{\rm rot}_\ep( {\bf \widetilde w}_\ep)\cdot {\rm rot}_\ep(\widetilde\psi_\ep)\,dx'dz_3+\ep^3 R_c\int_{\widetilde \Omega_\ep}{\rm div}_\ep({\bf \widetilde w}_\ep)\cdot {\rm div}_\ep(\widetilde\psi_\ep)\,dx'dz_3+4N^2\ep\int_{\widetilde \Omega_\ep}{\bf \widetilde w}_\ep\cdot\widetilde\psi_\ep\,dx'dz_3
\nonumber\\
\noame
&&\displaystyle -2N^2\ep\int_{\widetilde\Omega_\ep}{\bf\widetilde u}_\ep\cdot {\rm rot}_\ep(\widetilde\psi_\ep)\,dx'dz_3
-2N^2\beta\int_{\Gamma_0}({\bf \widetilde u_\ep}-{\bf s})\times {\bf n}\cdot\widetilde\psi_\ep\,d\sigma(x')\label{Form_Var_micro_tilde_hat}\\
\noame
&&\displaystyle =
-2N^2\int_{\Gamma_0}({\bf \widetilde u}_\ep\times{\bf n})\cdot\widetilde\psi_\ep\,d\sigma(x').\nonumber
\end{eqnarray}
Now, from the definition of ${\rm rot}_\ep$ given in (\ref{def_operator_tilde}), and by the change of variables given in Remark \ref{rem:hat} (see \cite{Anguiano_SG} for more details) applied to (\ref{Form_Var_vel_tilde_hat}), we obtain
\begin{equation}\label{Form_Var_vel_tilde_hat_v}
\begin{array}{l}
\displaystyle {\ep^2\over \ep^{2\ell}}\int_{\omega\times Z}{\rm rot}_{z'}(\widehat u_{3,\ep})\cdot {\rm rot}_{z'}(\widehat\varphi_{3})\,dx'dz+{\ep^2\over \ep^2}\int_{\omega\times Z}{\rm rot}_{z_3}({\bf\widehat  u}'_\ep)\cdot {\rm rot}_{z_3}(\widehat\varphi')\,dx'dz+{\ep^2\over \ep^{2\ell}}\int_{\omega\times Z}{\rm Rot}_{z'}({\bf\widehat  u}'_\ep)\cdot {\rm Rot}_{z'}(\widehat\varphi)\,dx'dz\\
\noame
\displaystyle-\ep^2\int_{\omega\times \Pi}\widehat P_\ep\,{\rm div}_{x'}(\widehat\varphi')\,dx'dz-{\ep^2\over \ep^\ell}\int_{\omega\times \Pi}\widehat P_\ep\,{\rm div}_{z'}(\widetilde\varphi')\,dx'dz-\ep\int_{\omega\times \Pi}\widehat P_\ep\,\partial_{z_3}\widehat\varphi\,dx'dz\\
\noame
\displaystyle -2N^2{\ep^2\over \ep^\ell}\int_{\omega\times Z} {\rm rot}_{z'}(\widehat\varphi_3)\cdot {\bf \widehat w}_\ep'\,dx'dz-2N^2\ep\int_{\omega\times Z} {\rm rot}_{z_3}(\widehat\varphi')\cdot {\bf \widehat w}'_\ep\,dx'dz-2N^2{\ep^2\over \ep^\ell}\int_{\omega\times Z} {\rm Rot}_{z'}(\widehat\varphi')\cdot \widetilde w_{3,\ep}\,dx'dz\\
\noame
\displaystyle +2\left(N^2-{1\over \alpha}\right)\ep\int_{\omega\times \widehat \Gamma_0}({\bf \widehat w}_\ep\times {\bf n})\cdot\widehat \varphi^b_\ep\,dx'd\sigma(z')+O_\ep=0,
\end{array}
\end{equation}
where $O_\ep$ is devoted to tends to zero. To simplify the variational formulation, we observe that from estimates for $({\bf \widehat u}_\ep, {\bf \widehat w}_\ep)$ given in Lemma \ref{estimates_hat} and convergence (\ref{RelationAlpha}), we deduce 
$$\left|{\ep^2\over \ep^{2\ell}}\int_{\widetilde \Omega_\ep}{\rm rot}_{z'}(\widehat u_{3,\ep})\cdot {\rm rot}_{z'}(\widehat\varphi_{3})\,dx'dz+{\ep^2\over \ep^{2\ell}}\int_{\widetilde \Omega_\ep}{\rm Rot}_{z'}({\bf\widehat  u}'_\ep)\cdot {\rm Rot}_{z'}(\widehat\varphi)\,dx'dz\right|\leq C\ep^{1-\ell}\to 0,$$
$$\left|2N^2{\ep^2\over \ep^\ell}\int_{\omega\times Z} {\rm rot}_{z'}(\widehat\varphi_3)\cdot {\bf \widehat w}_\ep'\,dx'dz+2N^2{\ep^2\over \ep^\ell}\int_{\omega\times Z} {\rm Rot}_{z'}(\widehat\varphi')\cdot \widetilde w_{3,\ep}\,dx'dz\right|\leq C\ep^{1-\ell}\to 0.$$
Then, variational formulation (\ref{Form_Var_vel_tilde_hat_v}) reads  
\begin{equation}\label{Form_Var_vel_tilde_hat_v2}
\begin{array}{l}
\displaystyle  \int_{\omega\times Z}{\rm rot}_{z_3}({\bf\widehat  u}'_\ep)\cdot {\rm rot}_{z_3}(\widehat\varphi')\,dx'dz-\ep^2\int_{\omega\times \Pi}\widehat P_\ep\,{\rm div}_{x'}(\widehat\varphi')\,dx'dz\\
\noame
\displaystyle-\ep^{2-\ell}\int_{\omega\times \Pi}\widehat P_\ep\,{\rm div}_{z'}(\widehat\varphi')\,dx'dz-\ep\int_{\omega\times \Pi}\widehat P_\ep\,\partial_{z_3}\widehat\varphi\,dx'dz\\
\noame
\displaystyle -2N^2\ep\int_{\omega\times Z} {\rm rot}_{z_3}(\widehat\varphi')\cdot {\bf \widehat w}'_\ep\,dx'dz +2\left(N^2-{1\over \alpha}\right)\ep\int_{\omega\times \widehat \Gamma_0}({\bf \widehat w}_\ep\times {\bf n})\cdot \widehat \varphi \,dx'd\sigma(z') +O_\ep=0,
\end{array}
\end{equation}
To finish, proceeding similarly for the variational formulation (\ref{Form_Var_micro_tilde_hat}), and using the estimates for $({\bf \widehat u}_\ep, {\bf \widehat w}_\ep)$ given in Lemma \ref{estimates_hat} and convergence (\ref{RelationAlpha}) as in previous variational formulation, we deduce that (\ref{Form_Var_micro_tilde_hat}) reads
\begin{equation}\label{Form_Var_vel_tilde_hat_w2}
\begin{array}{l}
\displaystyle
\ep R_c\int_{\omega\times Z}{\rm rot}_{z_3}( {\bf \widehat w}'_\ep)\cdot {\rm rot}_{z_3}(\widehat\psi')\,dx'dz+\ep R_c\int_{\omega\times Z}\partial_{z_3}\widehat w_{3,\ep}\, \partial_{z_3}\widehat\psi_3\,dx'dz+4N^2\ep\int_{\omega\times Z}{\bf \widehat w}_\ep\cdot\widehat\psi\,dx'dz\\
\noame
\displaystyle -2N^2\int_{\omega\times Z}{\bf\widehat u}'_\ep\cdot {\rm rot}_{z_3}(\widehat\psi')\,dx'dz
-2N^2\beta\int_{\omega\times\widehat \Gamma_0}({\bf \widehat u_\ep}-{\bf s})\times {\bf n}\cdot  \widehat\psi \,dx'd\sigma(z') \\
\noame
\displaystyle =
-2N^2\int_{\omega\times\widehat\Gamma_0}({\bf \widehat u}_\ep \times{\bf n})\cdot \widehat\psi \,dx'd\sigma(z')+O_\ep,
\end{array}
\end{equation}
where $O_\ep$ is devoted to tends to zero.  

\subsection{Convergences and limit problem}
We give a compactness result concerning the related unfolding functions $({\bf \widehat u}_\ep, {\bf \widehat w}_\ep, \widehat P_\ep)$ satisfying the {\it a priori} estimates given in Lemmas  \ref{estimates_hat}.
\begin{lemma}\label{lem_asymp_sub_hat} Consider the functions $({\bf \widetilde u}, {\bf \widetilde w}, \widetilde p)$ obtained in Lemma \ref{lem_asymp_sub}. Then, for a subsequence of $\ep$ still denote by $\ep$, we have the following convergence results:
\begin{itemize}
\item[(i)] {\it (Velocity)} There exists ${\bf \widehat u}=({\bf \widehat u}', \widehat u_{3})\in H^1(0,h(z');L^2_\#(\omega\times Z')^3)$, with $\widehat u_3\equiv 0$ and ${\bf \widehat u}'=0$ on $\widehat \Gamma_1$, such that it holds 
\begin{equation}\label{relation_hat_tilde}
\int_{Z}{\bf \widehat u}'(x',z)dz =\int_0^{h_{\rm max}}{\bf \widetilde u}'(x',z_3)\,dz_3\,,
\end{equation} and moreover
\begin{eqnarray}
 {\bf \widehat u}_\ep\rightharpoonup {\bf \widehat u}&&\hbox{  in  }H^1(0,h(z');L^2(\omega\times Z')^3),\label{conv_u_sub_hat}
 \\
\noame
\displaystyle {\rm div}_{z'}{\bf \widehat u}'=0&&\hbox{  in  }\omega\times Z,\label{div_sub_hat1a}
\\
\noame
\displaystyle {\rm div}_{z'}\left(\int_0^{h(z')}{ \widehat u}'\,dz_3\right)=0&&\hbox{  in  }\omega\times Z',\label{div_sub_hat1}\\
\noame
\displaystyle {\rm div}_{x'}\left(\int_{Z}\widehat u'(x',z)\,dz\right)=0&& \hbox{  in  }\omega\,,\label{div_sub_hat2}\\
\noame
\displaystyle
\left(\int_{Z}\widehat u'(x',z)\,dz\right)\cdot {\bf n}'=0&& \hbox{  on  }\partial\omega\,.\nonumber
\end{eqnarray}
\item[(ii)] {\it (Microrotation)} There exists ${\bf \widehat w}\in H^1(0,h(z');L^2_\#(\omega\times Z')^3)$, with $\widehat w_3\equiv 0$ and ${\bf \widehat w}'=0$ on $\widehat \Gamma_1$, such that it holds
\begin{equation}\label{relation_hat_tilde_w}
\int_{Z}{\bf \widehat w}'(x',z)dz=\int_0^{h_{\rm max}}{\bf \widetilde w}'(x',z_3)\,dz_3,
\end{equation} and moreover
\begin{eqnarray}
&\ep{\bf \widehat w}_\ep\rightharpoonup {\bf \widehat w}\quad\hbox{  in  }H^1(0,h(z');L^2(\omega\times Z')^3).&\label{conv_w_sub_hat}
\end{eqnarray}
\item[(iii)] {\it (Pressure)} The following convergence holds
\begin{eqnarray}
&\displaystyle  \epsilon^{2} \widehat P_\ep\to  \widetilde p\hbox{  in  }L^2(\omega\times \Pi).\label{conv_P_hat}
\end{eqnarray}
\end{itemize}
\end{lemma}
\begin{proof} 
The proof follows the lines of the proof of Lemma \ref{lem_asymp_sub}, so we will give some remarks. We start by $(i)$. Estimates given in (\ref{estim_u_hat}) imply the existence of $\hat u\in  H^1(0,h(z');L^2(\omega\times Z')^3)$ such that convergence (\ref{conv_u_sub_hat}) holds. As in Lemma \ref{lem_asymp_sub}, it holds that ${\bf \widehat u}=0$ on $\omega\times \widehat \Gamma_1$ and ${\widehat u}_3=0$ on $\omega\times \widehat \Gamma_0$. It can also be proved the $Z'$-periodicity of $\hat u$. This can be obtained by proceeding as in
\cite[Lemma 5.4]{grau1}.\\

 We have to take into account that applying the unfolded change of variables  to the divergence condition ${\rm div}_{\ep}({\widetilde u}_\epsilon)=0$ and multiplying by $\ep$, we get  
\begin{equation}\label{div_proof}
\epsilon^{1-\ell}{\rm div}_{z'}({\bf \widehat u}'_\ep) +\partial_{z_3}\widehat u_{3,\ep}=0.
\end{equation}
Passing to the limit, since relation (\ref{RelationAlpha}), we get $\partial_{z_3}\widehat u_{3}=0$, which means that $\widehat u_3$ is independent of $z_3$. Due to the boundary conditions ${\widehat u}_3=0$ on $\omega\times (\widehat \Gamma_1 \cup \widehat \Gamma_0)$, it holds that $\widehat u_3\equiv 0$.  
\\

Now, multiplying (\ref{div_proof})  by $\epsilon^{\ell-1}\widetilde q$ with $\widetilde q$ independent of $z_3$, after integrating by parts, we get
$$\int_{\omega\times Z} {\bf \widehat u}'_\ep\cdot \nabla_{z'}\widetilde q\,dx'dz=0.$$
Passing to the limit and integrating by parts, we get (\ref{div_sub_hat1a}). 
Observe that previous equality also can be written as follows
$$\int_{\omega\times Z'}\left(\int_0^{h(z')}{\bf \widehat u}'_\ep\,dz_3\right)\nabla_{z'}\widetilde q\,dx'dz'=0,$$
which passing to the limit and integrating by parts, gives (\ref{div_sub_hat1}). 
\\

For the proof of (\ref{relation_hat_tilde}),  we refer to \cite[Lemma 5.4]{Anguiano_SG}. Next, putting relation (\ref{relation_hat_tilde}) into the divergence condition (\ref{div_x_sub_tilde}), we get (\ref{div_sub_hat2}).\\

The proofs of $(ii)$ is similar to the proof of $(i)$ just taking into account the estimates of ${\bf \widehat w}_\ep$. \\

Finally, to prove $(iii)$, we remark that the strong convergence of sequence $\ep^2\widehat P_\ep$ to $\widetilde p$ is a consequence of the strong convergence of $\ep^2 \widetilde P_\ep$  to $\widetilde p$, see \cite[Proposition 2.9]{Ciora2}.

\end{proof}

\noindent Using   convergences given in Lemma \ref{lem_asymp_sub}, we give the reduced two-pressured homogenized system satisfied by $({\bf \widehat u},{\bf \widehat w}, \widetilde p)$.

\begin{theorem}[Limit problem]\label{thm_sub1}
The triplet  of functions 
   $({\bf \widehat u}, {\bf \widehat w}, \widetilde p)\in H^1(0,h(z');$ $L^2_\#(\omega\times Z')^3)\times H^1(0,h(z');L^2_\#(\omega\times Z')^3)\times  (L^2_0(\omega)\cap H^1(\omega))$, with $\widehat u_3=\widehat w_3\equiv 0$, given in Lemma \ref{lem_asymp_sub} is the unique solution of the two-pressure homogenized reduced micropolar problem
\begin{equation}\label{hom_system_sub_u}
\left\{\begin{array}{rcll}
\displaystyle
-\partial_{z_3}^2 {\bf \widehat u}'+ \nabla_{z'}\widehat \pi -2N^2{\rm rot}_{z_3}({\bf\widehat  w}')&=&\displaystyle - \nabla_{x'}\widetilde  p(x')&\hbox{ in }\omega\times Z,\\
\noame
\displaystyle
-R_c\partial_{z_3}^2 {\bf \widehat w}'+4N^2{\bf\widehat  w}' - 2N^2{\rm rot}_{z_3}({\bf \widehat u}')&=&0\displaystyle &\hbox{ in }\omega\times Z,\\
\noame
\widehat \pi \in L^2_{\#}(\omega\times Z'),
\end{array}\right.
\end{equation}
with divergence conditions
\begin{eqnarray} 
\displaystyle {\rm div}_{z'}({\bf \widehat u}')=0&&\hbox{ in }\omega\times Z,\label{divz_limit}\\
\noame
\displaystyle {\rm div}_{z'}\left(\int_0^{h(z')}{\bf \widehat u}'(x',z)\,dz_3\right)=0&&\hbox{ in }\omega\times Z',\label{divzh_limit}\\
\noame
\displaystyle {\rm div}_{x'}\left(\int_{Z}{\bf \widehat u}'(x_1,z)\,dz\right)=0&&\hbox{ in }\omega,\label{divx_limit}
\\
\noame\displaystyle
\left(\int_{Z}{\bf \widehat u}'(x',z)\,dz\right)\cdot {\bf n}'=0&& \hbox{  on  }\partial\omega\,,\nonumber
\end{eqnarray}

and boundary conditions
\begin{eqnarray} 
{\bf \widehat u}'=0,\quad{\bf \widehat w}'=0&&\hbox{on }\omega\times \widehat\Gamma_1,\label{bc_limit1}\\
\noame
\displaystyle \partial_{z_3}{\bf \widehat u}'=-{2\over \alpha}({\bf \widehat w}')^\perp&&\hbox{on }\omega\times \widehat\Gamma_0,\label{bc_limit2}\\
\noame
\displaystyle \partial_{z_3}{\bf \widehat w}'=-2N^2 \beta({\bf \widehat u}'-{\bf s})^\perp&&\hbox{on }\omega\times \widehat\Gamma_0.\label{bc_limit3}
\end{eqnarray}
\end{theorem}
\begin{proof} From Lemma \ref{lem_asymp_sub}, it remains to prove (\ref{hom_system_sub_u}) and boundary conditions (\ref{bc_limit2}) and (\ref{bc_limit3}). We divide the proof in three steps.  \\

\noindent {\it Step 1.} We prove  (\ref{hom_system_sub_u})$_{1}$ with boundary condition (\ref{bc_limit2}). According to Lemma \ref{lem_asymp_sub}, we consider  in (\ref{Form_Var_vel_tilde_hat_v2}) where $\widehat \varphi(x',z)\in \mathcal{D}(\omega;C^\infty_{\#}(Z)^2)$ with $\widehat \varphi_3\equiv 0$ in $\omega\times Z$ and ${\rm div}_{z'} (\widehat \varphi')=0$ in $\omega\times Z$. This gives
\begin{equation}\label{Form_Var_vel_tilde_hat_v2_proof1}
\begin{array}{l}
\displaystyle  \int_{\omega\times Z}{\rm rot}_{z_3}({\bf\widehat  u}'_\ep)\cdot {\rm rot}_{z_3}(\widehat\varphi')\,dx'dz-\ep^2\int_{\omega\times \Pi}\widehat P_\ep\,{\rm div}_{x'}(\widehat\varphi')\,dx'dz\\
\noame
\displaystyle -2N^2\ep\int_{\omega\times Z} {\rm rot}_{z_3}(\widehat\varphi')\cdot {\bf \widehat w}'_\ep\,dx'dz +2\left(N^2-{1\over \alpha}\right)\ep\int_{\omega\times \widehat \Gamma_0}({\bf \widehat w}_\ep\times {\bf n})\cdot  \widehat \varphi^b_\ep\,dx'd\sigma(z') +O_\ep=0,
\end{array}
\end{equation}
where $O_\epsilon$ is devoted to tends to zero when $\epsilon\to 0$. Below, let us pass to the limit when $\epsilon$ tends to zero in each term of the previous variational formulation:
\begin{itemize}
\item First term. Using convergence (\ref{conv_u_sub_hat}) and the definition of the operator ${\rm rot}_{z_3}$ given in (\ref{rotationals}), we get
$$\begin{array}{rl}
\displaystyle \int_{\omega\times Z}{\rm rot}_{z_3}({\bf\widehat  u}'_\ep)\cdot {\rm rot}_{z_3}(\widehat\varphi')\,dx'dz=&\displaystyle 
\int_{\omega\times Z}{\rm rot}_{z_3}({\bf\widehat  u}')\cdot {\rm rot}_{z_3}(\widehat\varphi')\,dx'dz+O_\ep\\
\noame
=& \displaystyle
\int_{\omega\times Z}\partial_{z_3}{\bf\widehat  u}'\cdot \partial_{z_3}\widehat\varphi'\,dx'dz+O_\ep.
\end{array}
$$
\item Second term. Using convergence (\ref{conv_P_hat}), we get
$$\ep^2\int_{\omega\times \Pi}\widehat P_\ep\,{\rm div}_{x'}(\widehat\varphi')\,dx'dz= \int_{\omega\times \Pi}\widetilde p\,{\rm div}_{x'}(\widehat\varphi')\,dx'dz+O_\ep.$$
\item Third term. Using convergence (\ref{conv_w_sub_hat}) and integration by parts, we get
$$\begin{array}{rl}\displaystyle
\displaystyle 
-2N^2\ep\int_{\omega\times Z} {\rm rot}_{z_3}(\widehat\varphi')\cdot {\bf \widehat w}'_\ep\,dx'dz 
=&\displaystyle 
-2N^2 \int_{\omega\times Z} {\rm rot}_{z_3}(\widehat\varphi')\cdot {\bf \widehat w}'\,dx'dz +O_\ep\\
\noame
=&\displaystyle 
-2N^2 \int_{\omega\times Z} {\rm rot}_{z_3}({\bf \widehat w}')\cdot \widehat\varphi'\,dx'dz+2N^2 \int_{\omega\times \widehat \Gamma_0}  (\widehat\varphi')^\perp\cdot {\bf \widehat w}'dx'd\sigma(z')+O_\ep\\
\noame
=&\displaystyle 
-2N^2 \int_{\omega\times Z} {\rm rot}_{z_3}({\bf \widehat w}')\cdot \widehat\varphi'\,dx'dz-2N^2 \int_{\omega\times \widehat \Gamma_0}  ({\bf \widehat w}')^\perp\cdot \widehat\varphi'dx'd\sigma(z')+O_\ep.
\end{array}$$
\item Fouth term. By continuity of the trace operator from $H^1(0,h(z');L^2(\omega\times Z')^3)$ into $L^2(\omega\times \widehat \Gamma_0)$ and using convergence  (\ref{conv_w_sub_hat}), we get the convergence of $\ep {\bf \widehat w}_\ep(x',z'0)$  to ${\bf \widehat w}$, and so we get
$$\begin{array}{rl}\displaystyle
2\left(N^2-{1\over \alpha}\right)\ep\int_{\omega\times \widehat \Gamma_0}({\bf \widehat w}_\ep\times {\bf n})\cdot  \widehat \varphi \,dx'd\sigma(z')=&\displaystyle2\left(N^2-{1\over \alpha}\right) \int_{\omega\times \widehat \Gamma_0}({\bf \widehat w}\times {\bf n})\cdot\widehat \varphi\,dx'd\sigma(z')+O_\ep\\
\noame
=&\displaystyle
2\left(N^2-{1\over \alpha}\right)\int_{\omega\times \widehat \Gamma_0}({\bf \widehat w}')^\perp\cdot\widehat \varphi'\,dx'd\sigma(z')+O_\ep.
\end{array}$$
\end{itemize}
Therefore, by previous convergences, we deduce that the limit variational formulation is given by the following one
\begin{equation}\label{form_var_1_changvar_hat1_limit2_sub2}
 \begin{array}{l}
\displaystyle \int_{\omega\times Z}\partial_{z_3}{\bf\widehat  u}'\cdot \partial_{z_3}\widehat\varphi'\,dx'dz+ \int_{\omega\times \Pi}\widetilde p\,{\rm div}_{x'}(\widehat\varphi')\,dx'dz\\
\noame
\displaystyle 
-2N^2 \int_{\omega\times Z} {\rm rot}_{z_3}({\bf \widehat w}')\cdot \widehat\varphi'\,dx'dz-{2\over \alpha}\int_{\omega\times \widehat \Gamma_0}({\bf \widehat w}')^\perp\cdot\widehat \varphi'\,dx'd\sigma(z')=0,
\end{array}\end{equation}
for every $\widehat \varphi'\in\mathcal{D}(\omega;C^\infty_\#(Z)^2)$ with ${\rm div}_{z'} (\widehat \varphi')=0$ in $\omega\times Z$.  If we consider that ${\rm div}_{x'}\left(\int_Z\widehat \varphi'(x',z)dz\right)=0$, then 
$$\int_{\omega\times \Pi}\widetilde p\,{\rm div}_{x'}(\widehat\varphi')\,dx'dz=\int_{\omega}\widetilde p\,{\rm div}_{x'}\left(\int_Z\widehat\varphi'\,dz\right)dx'=0,$$
and so, the previous variational formulation is  
\begin{equation}\label{form_var_1_changvar_hat1_limit2_sub2_111}
 \begin{array}{l}
\displaystyle \int_{\omega\times Z}\partial_{z_3}{\bf\widehat  u}'\cdot \partial_{z_3}\widehat\varphi'\,dx'dz
-2N^2 \int_{\omega\times Z} {\rm rot}_{z_3}({\bf \widehat w}')\cdot \widehat\varphi'\,dx'dz-{2\over \alpha}\int_{\omega\times \widehat \Gamma_0}({\bf \widehat w}')^\perp\cdot\widehat \varphi'\,dx'd\sigma(z')=0,
\end{array}\end{equation}
which, by density, holds for every function $\widehat \varphi' \in \mathcal{V}$ with 
$$\mathcal{V}_2=\left\{\begin{array}{l}
\displaystyle 
\displaystyle \widehat \varphi'\in H^1(0,h(z_1);L^2_\#(\omega\times Z')^2) :\qquad {\rm div}_{z'} ( \widehat \varphi')=0 \hbox{ in }\omega\times Z\\
\noame
\displaystyle {\rm div}_{x'}\left(\int_Z\widehat \varphi'(x',z)dz\right)=0\ \hbox{in }\omega,\quad \left(\int_Z\widehat \varphi'(x',z)dz\right)\cdot {\bf n}'=0\ \hbox{ on }\partial\omega\end{array}\right\}.$$
 Reasoning as in \cite{Allaire, Anguiano_evolutivo, Anguiano_SG}, the orthogonal of $\mathcal{V}_2$ with respect to the usual scalar product in $L^2(\omega \times Z')$ is made of gradients of the form $\nabla_{x'}\widetilde q(x')+\nabla_{z'}\widehat \pi(x',z')$ with $\widetilde q(x')\in L^2_0(\omega)$  and $\widehat \pi(x',z')\in L^2_\#(\omega\times Z')$.   Therefore, by integration by parts, the variational formulation (\ref{form_var_1_changvar_hat1_limit2_sub2_111}) is equivalent to problem  (\ref{hom_system_sub_u})$_{1}$ with boundary condition (\ref{bc_limit2}).  It remains to prove that the pressure $\widetilde p(x')$, arising as a Lagrange multiplier of the incompressibility condition ${\rm div}_{x'}(\int_Z \widehat u'(x',z)dz)=0$, is the same as the limit of the pressure $\widetilde p_\epsilon$. This can be easily done by considering a test function only with ${\rm div}_{z'}$ equal to zero, arriving to the variational formulation (\ref{form_var_1_changvar_hat1_limit2_sub2}) and identifying limits.\\

\noindent {\it Step 2.} We prove   (\ref{hom_system_sub_u})$_{2}$ with boundary condition (\ref{bc_limit3}). According to Lemma \ref{lem_asymp_sub}, we consider  in (\ref{Form_Var_vel_tilde_hat_v2}) where $\widehat \psi(x',z)\in \mathcal{D}(\omega;C^\infty_{\#}(Z)^2)$ with $\widehat \psi_3\equiv 0$ in $\omega\times Z$. This gives
\begin{equation}\label{Form_Var_vel_tilde_hat_w2_proof1}
\begin{array}{l}
\displaystyle
\ep R_c\int_{\omega\times Z}{\rm rot}_{z_3}( {\bf \widehat w}'_\ep)\cdot {\rm rot}_{z_3}(\widehat\psi')\,dx'dz+4N^2\ep\int_{\omega\times Z}{\bf \widehat w}_\ep'\cdot\widehat\psi'\,dx'dz\\
\noame
\displaystyle -2N^2\int_{\omega\times Z}{\bf\widehat u}'_\ep\cdot {\rm rot}_{z_3}(\widehat\psi')\,dx'dz
-2N^2\beta\int_{\omega\times\widehat \Gamma_0}({\bf \widehat u_\ep}-{\bf s})\times {\bf n}\cdot  \widehat\psi \,dx'd\sigma(z') \\
\noame
\displaystyle =
-2N^2\int_{\omega\times\widehat\Gamma_0}({\bf \widehat u}_\ep \times{\bf n})\cdot  \widehat\psi \,dx'd\sigma(z')+O_\ep,
\end{array}
\end{equation}
where $O_\epsilon$ is devoted to tends to zero when $\epsilon\to 0$. Below, let us pass to the limit when $\epsilon$ tends to zero in each term of the previous variational formulation:
\begin{itemize}
\item First term. Using convergence (\ref{conv_w_sub_hat}) and the definition of the operator ${\rm rot}_{z_3}$ given in (\ref{rotationals}), we get
$$\begin{array}{rl}
\displaystyle \ep R_c\int_{\omega\times Z}{\rm rot}_{z_3}( {\bf \widehat w}'_\ep)\cdot {\rm rot}_{z_3}(\widehat\psi')\,dx'dz=&\displaystyle 
 R_c\int_{\omega\times Z}{\rm rot}_{z_3}( {\bf \widehat w}')\cdot {\rm rot}_{z_3}(\widehat\psi')\,dx'dz+O_\ep\\
 \noame
= &\displaystyle
 R_c\int_{\omega\times Z}\partial_{z_3}{\bf \widehat w}'\cdot \partial_{z_3}\widehat\psi'\,dx'dz+O_\ep.
\end{array}$$
\item Second term. Using convergence (\ref{conv_w_sub_hat}), we get
$$\begin{array}{rl}
\displaystyle 4N^2\ep\int_{\omega\times Z}{\bf \widehat w}_\ep'\cdot\widehat\psi'\,dx'dz=&\displaystyle 
 4N^2\int_{\omega\times Z}{\bf \widehat w}'\cdot\widehat\psi'\,dx'dz+O_\ep.
\end{array}$$
\item Third term. Using convergence (\ref{conv_u_sub_hat}) and integration by parts, we get
$$\begin{array}{rl}
\displaystyle 
-2N^2\int_{\omega\times Z}{\bf\widehat u}'_\ep\cdot {\rm rot}_{z_3}(\widehat\psi')\,dx'dz=&\displaystyle 
-2N^2\int_{\omega\times Z}{\bf\widehat u}'\cdot {\rm rot}_{z_3}(\widehat\psi')\,dx'dz+O_\ep\\
\noame
=&\displaystyle -2N^2\int_{\omega\times Z}\widehat\psi'\cdot {\rm rot}_{z_3}({\bf\widehat u}')\,dx'dz+2N^2\int_{\omega\times \widehat \Gamma_0}(\widehat\psi')^\perp\cdot {\bf\widehat u}'\,dx'd\sigma(z')+      O_\ep\\
\noame
=&\displaystyle -2N^2\int_{\omega\times Z}\widehat\psi'\cdot {\rm rot}_{z_3}({\bf\widehat u}')\,dx'dz-2N^2\int_{\omega\times \widehat \Gamma_0}({\bf\widehat u}')^\perp\cdot \widehat\psi'\,dx'd\sigma(z')+      O_\ep.
\end{array}$$
\item Fourth term. By continuity of the trace operator from $H^1(0,h(z');L^2(\omega\times Z')^3)$ into $L^2(\omega\times \widehat \Gamma_0)$ and from convergence  (\ref{conv_u_sub_hat}), we get the convergence of ${\bf \widehat u}_\ep(x',z',0)$ to ${\bf \widehat u}(x',z',0)$ and so
$$\begin{array}{rl}
\displaystyle -2N^2\beta\int_{\omega\times\widehat \Gamma_0}({\bf \widehat u_\ep}'-{\bf s})\times {\bf n}\cdot  \widehat\psi \,dx'd\sigma(z')=&\displaystyle-2N^2\beta\int_{\omega\times\widehat \Gamma_0}({\bf \widehat u}-{\bf s})\times {\bf n}\cdot \widehat\psi\,dx'd\sigma(z')+O_\ep\\
\noame
=&\displaystyle-2N^2\beta\int_{\omega\times\widehat \Gamma_0}({\bf \widehat u}'-{\bf s}')^\perp\cdot \widehat\psi'\,dx'd\sigma(z')+O_\ep.
\end{array}
$$
\item Five term.  Similar to the previous term, from convergence  (\ref{conv_u_sub_hat}), we get
$$\begin{array}{rl}
\displaystyle-2N^2\int_{\omega\times\widehat\Gamma_0}({\bf \widehat u}_\ep \times{\bf n})\cdot  \widehat\psi \,dx'd\sigma(z')=&\displaystyle-2N^2\int_{\omega\times\widehat\Gamma_0}({\bf \widehat u}')^\perp \cdot \widehat\psi'\,dx'd\sigma(z')+O_\ep.
\end{array}$$
\end{itemize}
Therefore, by previous convergences, we deduce that the limit variational formulation is given by the following one
\begin{equation}\label{form_var_1_changvar_hat1_limit2_sub2_w}
 \begin{array}{l}
\displaystyle  R_c\int_{\omega\times Z}\partial_{z_3}{\bf \widehat w}'\cdot \partial_{z_3}\widehat\psi'\,dx'dz+4N^2\int_{\omega\times Z}{\bf \widehat w}'\cdot\widehat\psi'\,dx'dz\\
\noame
\displaystyle 
-2N^2\int_{\omega\times Z}\widehat\psi'\cdot {\rm rot}_{z_3}({\bf\widehat u}')\,dx'dz-2N^2\beta\int_{\omega\times\widehat \Gamma_0}({\bf \widehat u}'-{\bf s})^\perp\cdot \widehat\psi'\,dx'd\sigma(z')=0,
\end{array}\end{equation}
for every $\widehat \psi'\in\mathcal{D}(\omega;C^\infty_\#(Z)^2)$. By density, (\ref{form_var_1_changvar_hat1_limit2_sub2_w}) holds for every function $\psi'$ in the $H^1(0,h(z_1);L^2_\#(\omega\times Z'))$, and is equivalent to problem  (\ref{hom_system_sub_u})$_{2}$ with boundary condition (\ref{bc_limit3}).\\

{\it Step 3. Conclusion. }  Since $\widehat\varphi'$ and $\widehat\psi$ are arbitrary, we derive from  (\ref{form_var_1_changvar_hat1_limit2_sub2}) and (\ref{form_var_1_changvar_hat1_limit2_sub2_w}) that $({\bf \widehat u}', {\bf \widehat w}', \widetilde p, \widehat \pi)$ satisfies the system  (\ref{hom_system_sub_u}) with boundary conditions (\ref{bc_limit1})--(\ref{bc_limit3}). 
 To ensure that the whole sequence $({\bf \widehat u}_\ep, \ep {\bf \widehat w}_\ep, \ep^2\widehat P_\ep)$ converges, it remains to prove the existence and uniqueness of weak solution of the effective system  (\ref{hom_system_sub_u}). This follows the lines of the proof of Theorem \ref{Existence_condition}, so we omit it.
 
\end{proof}

\section{Reynolds equation}\label{sec:reynolds}
In this section, we give the main result of this paper, i.e. the Reynolds equation for the pressure $\widetilde p$ and the expressions for the average velocity and microrotation. We will proceed as follows. First, we give the expressions of ${\bf \widehat u}$ and ${\bf \widehat w}$ by solving problem (\ref{hom_system_sub_u}) with boundary conditions (\ref{bc_limit1})--(\ref{bc_limit3}). Next, we integrate in $Z$ the expression of ${\bf \widehat u}$ and ${\bf \widehat w}$, which gives the expression for the average velocity and microrotation, which are given depending on local problems. Finally, putting the expression of the average velocity into the incompressibility condition (\ref{divx_limit}), we deduce the Reynolds equation.

\begin{lemma}\label{expression_uw_hat_alpha}
Assuming condition (\ref{Existence_condition}) and $\alpha\neq 1$, the solutions of (\ref{hom_system_sub_u}) with boundary conditions (\ref{bc_limit1})--(\ref{bc_limit3}) are given by the following expressions
\begin{eqnarray}
&{\bf \widehat u}'(x',z)=&\displaystyle \left[{2N^2\over k}\Big(\sinh(kz_3)-\sinh(kh(z'))+\gamma_\alpha (z_3-h(z'))\Big)A_1(z') \right.\nonumber\\
\noame
&&\displaystyle  \left.+{2N^2\over k}\left(\cosh(kz_3)-\cosh(kh(z'))B_1(z')+{z_3^2-h(z')^2\over 2(1-N^2)}\right)\right]\left(\nabla_{x'}\widetilde p(x')
+\nabla_{z'}\widehat \pi(x',z')\right)\label{expression_hat_u}\\
\noame
&&\displaystyle+ \left[\Big({2N^2\over k}\left(\sinh(kz_3)-\sinh(kh(z'))+\gamma_\alpha (z_3-h(z')\right)\Big)A_2(z')\right.\nonumber\\
\noame
&&\displaystyle \left. +{2N^2\over k}\left(\cosh(kz_3)-\cosh(kh(z'))\right)B_2(z')\right]{\bf s}',\nonumber
\end{eqnarray}

\begin{eqnarray}
&{\bf \widehat w}'(x',z)=&\displaystyle \Big[\Big(\cosh(kz_3)+{\gamma_\alpha\over 2}\Big)A_1(z')+\sinh(kz_3)B_1(z')+{z_3\over 2(1-N^2)}\Big](\nabla_{x'}\widetilde p(x')+\nabla_{z'}\widehat \pi(x',z'))^\perp
\nonumber\\
\noame
&&\displaystyle  +\Big[\Big(\cosh(kz_3)-{\gamma_\alpha\over 2}\Big)A_2(z')+\sinh(kz_3)B_2(z')\Big]({\bf s}')^\perp,\label{expression_hat_w}
\end{eqnarray}
where 
\begin{equation}\label{parameterk}
k=2N\sqrt{{1-N^2\over R_c}},\quad\quad {\gamma_\alpha\over 2}={1-\alpha N^2\over \alpha-1},
\end{equation}
and  $A_i, B_2$, $i=1,2$, are given by 
\begin{eqnarray}
&A_1(z')=&{L(z')\over 2(1-N^2)}\Big[h(z')\Big(4N^4(1-\cosh(kh(z')))+{R_c\over \beta} k^2\Big)\nonumber\\
\noame
&&-k\sinh(kh(z'))\Big({R_c\over \beta}-2N^2h(z')^2\Big)\Big],\label{A1}
\\
\noame
&A_2(z')=&-2N^2k L(z')\sinh(kh(z')),\label{A2}
\\
\noame
&B_1(z')=&{L(z')\over 2(1-N^2)}\Big[2N^2 h(z')\Big(2N^2\sinh(kh(z'))+\gamma_\alpha kh(z')\Big)\nonumber\\
\noame
&&+k({R_c\over \beta}-2N^2 h(z')^2)\Big(\cosh(kh(z'))+{\gamma_\alpha\over 2}\Big)\Big],\label{B1}
\\
\noame
&B_2(z')=&-2N^2kL(z')\Big[{\gamma_\alpha\over 2}+\cosh(kh(z'))\Big],\label{B2}
\end{eqnarray}
and
\begin{equation}\label{L}
\begin{array}{rl}
L(z')=&-\Big[\Big({\gamma_\alpha\over 2}+\cosh(kh(z'))\Big)\Big(4N^4[1-\cosh(kh(z'))]+{R_c\over \beta} k^2\Big)\\
\noame
&+2N^2\sinh(kh(z'))\Big(\gamma_\alpha kh(z')+2N^2\sinh(kh(z'))\Big)\Big]^{-1}.
\end{array}\end{equation}
\end{lemma}
\begin{proof}
The proof is obtaining by solving system (\ref{hom_system_sub_u})  for $(\widehat u_1, \widehat w_2)$ with corresponding boundary conditions (\ref{bc_limit1})--(\ref{bc_limit3}), by observing that problem (\ref{hom_system_sub_u})  is an ordinary differential equation with respect to $z_3$ (i.e. considering variables $x',z'$ as parameters). The proof can be found in \cite[Lemma 3.5]{Bayada_NewModel} for a system (\ref{hom_system_sub_u}) defined only for variables $(x',z_3)$, but the proof is the same just taking into account that here, two pressures appear. Following the proof, we observe that variable $x'$ only appears in the pressures $\widetilde p(x')$ and $\widehat \pi(x',z')$. 

We note that  $(\widehat u_1, \widehat w_2)$ and $(\widehat u_2, -\widehat w_1)$  satisfy the same equations and boundary conditions by writting $s_2$ and $z_2$ instead of $s_1$ and $z_1$. So, expression for $(\widehat u_2, \widehat w_1)$ are obtained straightforward, from expression of $(\widehat u_1, \widehat w_2)$, see \cite[Appendix 1]{Bayada_NewModel2}.

\end{proof}

\begin{lemma}\label{expression_uw_hat_alpha1}
Assuming condition (\ref{Existence_condition}) and $\alpha= 1$, the solutions of (\ref{hom_system_sub_u}) with boundary conditions (\ref{bc_limit1})--(\ref{bc_limit3}) are given by the following expressions
\begin{eqnarray}
&{\bf \widehat u}'(x',z)=&\displaystyle \Big[{2N^2\over k}(\cosh(kz_3)-\cosh(kh(z'))B_1'(z')+{z_3^2-h^2(z')\over 2(1-N^2)}+{A_1'(z')\over 1-N^2}(z_3-h(z'))\Big](\nabla_{x'}\widetilde p(x')+\nabla_{z'}\widehat \pi(x',z')) \nonumber\\
\noame
&&\displaystyle  +\Big[{2N^2\over k}(\cosh(kz_3)-\cosh(kh(z'))B_2'(z')+{A_2'(z')\over 1-N^2}(z_3-h(z'))\Big]{\bf s}',\label{expression_hat_u_alpha1}\\
\noame
&{\bf \widehat w}(x',z)=&\displaystyle \Big[\sinh(kz_3)B_1'(z')+{z_3\over 2(1-N^2)}+{A_1'(z')\over 2(1-N^2)}\Big](\nabla_{x'}\widetilde p(x')+\nabla_{z'}\widehat \pi(x', z'))^\perp
\label{expression_hat_w_alpha1}\\
\noame
&&\displaystyle + \Big[\sinh(kz_3)B_2'(z')+{A_2'(z')\over 2(1-N^2)}\Big]({\bf s}')^\perp,\nonumber
\end{eqnarray}
where $k$ is given by (\ref{parameterk}) and  $A'_i, B'_i$, $i=1,2$, are given by 
\begin{eqnarray}
&A_1'(z')=&L'(z')\Big[h(z')\Big(4N^4(1-\cosh(kh(z')))+{R_c\over \beta}k^2\Big)-k\sinh(kh(z'))\Big({R_c\over \beta}-2N^2h^2(z')\Big)\Big],
\label{A11}
\\
\noame
&A_2'(z')=&4N^2k(1-N^2)L'(z')\sinh(kh(z')),
\label{A21}
\\
\noame
&B_1'(z')=&k{L(z')\over 2(1-N^2)}\Big[2N^2h(z')+{R_c\over \beta}\Big],
\label{B11}
\\
\noame
&B_2'(z')=&-2N^2kL'(z'),
\label{B21}
\end{eqnarray}
and
\begin{equation}\label{L1}
L(z')=-\Big[4N^4(1-\cosh(kh(z'))+{R_c\over \beta}k^2+4N^2k h(z')\sinh(kh(z'))\Big]^{-1}.
\end{equation}
\end{lemma}
\begin{proof} The proof is similar to the case $\alpha\neq 1$, see \cite[Lemma 3.6]{Bayada_NewModel} and \cite[Appendix 2]{Bayada_NewModel2}.

\end{proof}

\begin{lemma}\label{lem:av_z3}  The $z_3$-average  velocity and microrotation satisfy the following expressions:
\begin{eqnarray}
&\displaystyle{\bf \widehat U}'(x',z')= \int_{0}^{h(z')}{\bf \widehat u}'(x',z)\,dz_3=\displaystyle-\Theta_1(z')(\nabla_{x'}\widetilde p(x')+\nabla_{z'}\widehat \pi(z'))+\Theta_2(z'){\bf s}',&\widehat U_{3}\equiv 0,\label{av_vel}\\
\noame
& \displaystyle{\bf \widehat W}'(x',z')=  \int_{0}^{h(z')}{\bf \widehat w}'(x',z)\,dz_3=\displaystyle\Phi_1(z')(\nabla_{x'}\widetilde p(x')+\nabla_{z'}\widehat \pi(z'))^\perp+\Phi_2(z')({\bf s}')^\perp,&\widehat W_{3}\equiv 0.\label{av_mic}
\end{eqnarray}
Here, functions $\Theta_i$ and $\Phi_i$, $i=1,2$, are given depending on the case:
\begin{itemize}
\item[--] If  $\alpha\neq 1$, then
\begin{eqnarray}
& \Theta_1(z')=&  {h^3(z')\over 3(1-N^2)}\label{Theta1}\\
\noame
&&  -\left[ {2N^2\over k}\Big({\cosh(kh(z'))-1\over k}-h(z')\sinh(kh(z'))\Big)-{\gamma_\alpha\over 2}h^2(z')\right]A_1(z')\nonumber\\
\noame
&& -{2N^2\over k}\Big({\sinh(kh(z'))\over k}-h(z')\cosh(kh(z'))\Big)B_1(z'),\nonumber\\
\noame
& \Theta_2(z')=& \left[{2N^2\over k}\left({\cosh(kh(z'))-1\over k}-h(z')\sinh(kh(z'))\right)-{\gamma_\alpha\over 2}h^2(z')\right]A_2(z') \label{Theta2}\\
\noame
&&+{2N^2\over k}\left({\sinh(kh(z'))\over k}-h(z')\cosh(kh(z'))\right)B_2(z'),\nonumber\\
\noame
&\Phi_1(z')=& {h(z')^2\over 4(1-N^2)}+ \Big({\sinh(kh(z'))\over k}+{\gamma_\alpha\over 2}\Big)A_1(z')+{\cosh(kh(z'))-1\over k}B_1(z'),\label{Phi1}\\
\noame
&\Phi_2(z')=& \Big({\sinh(kh(z'))\over k}-{\gamma_\alpha h(z')\over 2}\Big)A_1(z')+{\cosh(kh(z'))-1\over k}B_2(z'),\label{Phi2}
\end{eqnarray}
where $A_i, B_i$, $i=1,2$ are defined in Lemma \ref{expression_uw_hat_alpha}.

\item[--] If  $\alpha= 1$, then
\begin{eqnarray}
& \Theta_1(z')=& {h^3(z')\over 3(1-N^2)}+{h^2(z')\over 2(1-N^2)}A'_1(z')  \label{Theta11}\\
\noame
&&-{2N^2\over k}\left[{\sinh(kh(z'))\over k}-h(z')\cosh(kh(z'))\right]B_1'(z'),\nonumber\\
\noame
& \Theta_2(z')=& -{h^2(z')\over 2(1-N^2)}A'_2(z') \label{Theta21}\\
\noame
&&+{2N^2\over k}\left[{\sinh(kh(z'))\over k}-h(z')\cosh(kh(z'))\right]B_2'(z'),\nonumber\\
\noame
&\Phi_1(z')=&{h^2(z')\over 4(1-N^2)}+{h(z')\over 2(1-N^2)}A_1'(z') +{\cosh(kh(z'))-1\over k}B_1'(z'),\label{Phi11}\\
\noame
&\Phi_2(z')=& {h(z')\over 2(1-N^2)}A_2'(z')+{\cosh(kh(z'))-1\over k}B_2'(z'),\label{Phi21}
\end{eqnarray}
where $A_i', B_i'$, $i=1,2$ are defined in Lemma \ref{expression_uw_hat_alpha1}.
\end{itemize}
\end{lemma}
\begin{proof} Using divergence condition (\ref{divzh_limit}), we obtain for a.e. $x'\in \omega$ that
$$\int_{ Z'}\left(\int_0^{h(z')}{\bf \widehat u}'(x',z)\,dz_3\right)\nabla_{z'}\theta(z')\,dx'dz'=0,\qquad\forall\,\theta \in H^1( Z').$$
From Lemmas \ref{expression_uw_hat_alpha} and \ref{expression_uw_hat_alpha1}, by averaging (\ref{expression_hat_u}) or (\ref{expression_hat_u_alpha1}) with respect to $z_3$ between $0$ and $h(z')$, we obtain
$${\bf \widehat U}'(x',z')=\int_0^{h(z')}{\bf \widehat u}'(x',z)\,dz_3=-\Theta_1(z')(\nabla_{x'}\widetilde p(x')+\nabla_{z'}\widehat \pi(z'))+{\bf s}'\,\Theta_2(z'),$$
i.e. (\ref{av_vel}) with $\Theta_i$, $i=1,2$ given by (\ref{Theta1}), (\ref{Theta2}) if $\alpha\neq 1$, and  by (\ref{Theta11}), (\ref{Theta21}) if $\alpha=1$. Similarly, we deduce (\ref{av_mic}) by averaging (\ref{expression_hat_w}) or (\ref{expression_hat_w_alpha1}) with respect to $z_3$ between $0$ and $h(z')$ as follows
$${\bf \widehat W}'(x',z')=\int_0^{h(z')}{\bf \widehat w}'(x',z)\,dz_3=\Phi_1(z')(\nabla_{x'}\widetilde p(x')+\nabla_{z'}\widehat \pi(z'))^\perp+({\bf s}')^\perp\Phi_2(z'),$$
with $\Phi_i$, $i=1,2$ given by (\ref{Phi1})--(\ref{Phi2}) if $\alpha\neq 1$, and  by (\ref{Phi11})--(\ref{Phi21}) if $\alpha=1$.

\end{proof}
Finally, we give the main result of this paper, i.e. the derivation of a Reynolds equation satisfied by $\widetilde p$ with the effects of the roughness of the top boundary and the non-standard boundary condition on the flat boundary of the original thin domain.
\begin{theorem}(Main result)\label{thm_main} Assuming condition (\ref{Existence_condition}) and considering the functions $({\bf \widehat u}, {\bf \widehat w})$ given in Lemmas \ref{expression_uw_hat_alpha} and \ref{expression_uw_hat_alpha1}, we have that the $z$-average velocity and microrotation 
$$  {\bf \widetilde U}(x')=\int_{Z'}{\bf \widehat U}(x',z')\,dz',\quad {\bf \widetilde W}(x')=\int_{Z'}{\bf \widehat W}(x',z')\,dz',$$
satisfy the following expressions
\begin{equation}\label{AverageVelMicro}
\begin{array}{rl}
{\bf \widetilde U}(x')=&\displaystyle K^{(1)}\nabla_{x'}\widetilde p(x')+L^{(1)}{\bf s}',\quad \widetilde U_{3}\equiv 0,\\
\\
{\bf \widetilde W}(x')=&\displaystyle K^{(2)}\nabla_{x'}\widetilde p(x')^\perp+L^{(2)}({\bf s}')^\perp,\quad \widetilde W_{3}\equiv 0,
\end{array}
\end{equation}
for a.e. $x'\in\omega$. Here,  $K^{(k)}\in \mathbb{R}^{2\times 2}$ and $L^{(k)}\in\mathbb{R}$, $k=1,2$, are defined by
\begin{equation}\label{KL1}
\begin{array}{rl}\displaystyle 
K^{(1)}=\int_{Z'}\Theta_1(z')\left(\begin{array}{cc}
\partial_{z_1}q^1(z')+1 & \partial_{z_1}q^2(z')\\
\noame
\partial_{z_2}q^1(z') & \partial_{z_2}q^2(z')+1
\end{array}\right)dz',& \displaystyle L^{(1)}=\int_{Z'}\Theta_2(z')\,dz',\\
\\
\displaystyle K^{(2)}=\int_{Z'}\Phi_1(z')\left(\begin{array}{cc}
\partial_{z_2}q^2(z')+1& -\partial_{z_2}q^1(z')\\
\noame
- \partial_{z_1}q^2(z')& \partial_{z_1}q^1(z')+1 
\end{array}\right)dz',&
\displaystyle L^{(2)}=\int_{Z'}\Phi_2(z')\,dz',
\end{array}
\end{equation}
where  function $q^i\in H^1_\#(Z')$, $i=1,2$, satisfies the   local problem 
\begin{equation}\label{LocalReynolds}
\int_{Z'}\Theta_1(z')(\nabla_{z'}q^i(z')+{\bf e}_i)\cdot \nabla_{z'}\theta(z')\,dz'=\int_{Z'}\Theta_2(z')s_i\,({\bf e}_i\cdot\nabla_{z'}\theta(z'))\,dz'\quad \forall\,\theta\in H^1(Z').
\end{equation}
Moreover, $\widetilde p$ satisfies the following Reynolds equation
\begin{equation}\label{Reynolds_problem}
 \begin{array}{rl}\displaystyle 
\int_{\omega}K^{(1)}\nabla_{x'}\widetilde p(x')\cdot \nabla_{x'}\eta(x')\,dx'=\int_{\omega}L^{(1)}{\bf s}'\cdot \nabla_{x'}\eta(x')\,dx'\quad \forall\,\eta\in H^1(\omega).
\end{array} 
\end{equation}
\end{theorem}
\begin{proof} From the expressions ${\bf \widehat U}'$ in (\ref{av_vel}) given by 
$${\bf \widehat U}'(x',z')=-\Theta_1(z')(\nabla_{x'}\widetilde p(x')+\nabla_{z'}\widehat \pi(z'))+\Theta_2(z'){\bf s}',$$
and divergence condition (\ref{divzh_limit}), i.e. ${\rm div}_{z'}({\bf \widehat U}')=0$ in $\omega\times Z'$, we know that $\nabla_{z'}\widehat \pi$ holds the following problem  
\begin{equation}\label{var_ppi}
\int_{\omega\times Z'}\Theta_1(z')(\nabla_{x'}\widetilde p(x')+\nabla_{z'}\widehat \pi(z'))\cdot \nabla_{z'}\vartheta(x',z')\,dx'dz'=\int_{\omega\times Z'}\Theta_2(z'){\bf s}'\cdot \nabla_{z'}\vartheta(x',z')\,dx'dz'\quad \forall \vartheta\in H^1(\omega\times Z').
\end{equation}
Now, to find a problem for $\nabla_{x'}\widetilde p$ we proceed to eliminate the microscopic variable $z'$. To do this, we define
$$\widehat \pi(x',z')=\sum_{i=1}^2\partial_{x_i}\widetilde p(x') q^i(z')\quad \hbox{in }\omega\times Z',$$
and so, it holds
$$\nabla_{x'}\widetilde p(x')+\nabla_{z'}\widehat \pi(z')=\left(\begin{array}{cc}
\partial_{z_1}q^1(z')+1 & \partial_{z_1}q^2(z')\\
\noame
\partial_{z_2}q^1(z') & \partial_{z_2}q^2(z')+1
\end{array}\right)\nabla_{x'}\widetilde p(x'),$$
and observe that  this choice of $\widehat \pi$ and taking $\vartheta(x',z')=\eta(x')\theta(z')$ with $\eta\in H^1(\omega)$ and $\theta\in H^1(Z')$, then $\widehat \pi$  satisfies (\ref{var_ppi}) by taking into account local problems (\ref{LocalReynolds}) for $q^i$, $i=1,2$. Then, integrating ${\bf \widehat U}_{\rm av}(x',z')$ in $Z'$, we deduce expression for ${\bf \widetilde U}(x')=\int_{Z'}{\bf \widehat U}(x',z')\,dz'$  given in (\ref{AverageVelMicro}).  

Finally, putting expression of ${\bf \widetilde U}$ in the divergence condition (\ref{divx_limit}), i.e. ${\rm div}_{x'}({\bf \widetilde U}')=0$ in $\omega$, we deduce that $\widetilde p$ satisfies the Reynolds equation (\ref{Reynolds_problem}).\\

The derivation of the expression of ${\bf \widetilde W}'$ given in (\ref{AverageVelMicro})$_2$ is straightforward   taking into account (\ref{av_mic}) and that 
$$(\nabla_{x'}\widetilde p(x')+\nabla_{z'}\widehat \pi(z'))^\perp=\left(\begin{array}{cc}
\partial_{z_2}q^2(z')+1& -\partial_{z_2}q^1(z')\\
\noame
- \partial_{z_1}q^2(z')& \partial_{z_1}q^1(z')+1 
\end{array}\right)(\nabla_{x'}\widetilde p(x'))^\perp.$$
Then, integrating ${\bf \widehat W}'(x',z')$ in $Z'$, we deduce expression for ${\bf \widetilde W}(x')=\int_{Z'}{\bf \widehat W}(x',z')\,dz'$  given in (\ref{AverageVelMicro})$_2$.

\end{proof}

\begin{remark} We make an observation regarding the case in which the non-standard boundary conditions are assumed on the rough surface $\Gamma_1^\ep=\omega\times \{h_\ep(x')\}$ and not on the flat one $\Gamma_0$.  In this case,  we need to obtain a relation between $\|{\rm rot}(\varphi)\|_{L^2(\Omega^\ep)^3}$ and $\|D\varphi\|_{L^2(\Omega^\ep)^3}$. We recall that it holds
\begin{equation}\label{remark_ref2}\int_{\Omega^\ep}(|{\rm div}(\varphi)|^2+|{\rm rot}(\varphi)|^2)dx=\int_{\Omega^\ep}|D\varphi|^2dx +\int_{\Gamma_1^\ep}((\varphi\cdot \nabla){\bf n}_\ep)\cdot \varphi\,d\sigma,
\end{equation}
for every $\varphi\in {\bf V}_\ep=\{\varphi\in H^1(\Omega^\ep)^3\,:\, \varphi_{\Gamma_0\cup \Gamma_{\rm lat}^\ep}=0, \varphi\cdot {\bf n}_\ep=0\ \hbox{ on }\Gamma_1^\ep\}$  (see, for instance \cite[formula (IV.23)]{Boyer}). In particular, if the boundary $\Gamma_1^\ep$ were flat, the identity $\|{\rm rot}(\varphi)\|_{L^2(\Omega^\ep)^3}=\|D\varphi\|_{L^2(\Omega^\ep)^3}$ would hold for every $\varphi\in {\bf V}_0$,  with ${\bf V}_0=\{\varphi\in {\bf V}_\ep\,:\, {\rm div}(\varphi)=0\ \hbox{in }\Omega^\ep\}$, since the remaining therm $\int_{\Gamma_1^\ep}((\varphi\cdot \nabla){\bf n}_\ep)\cdot \varphi\,d\sigma$ would vanish. However, in the present geometric configuration, one can not expect this term to be zero in general. In fact, a classical estimate reads
$$\left|\int_{\Gamma_1^\ep}((\varphi\cdot \nabla){\bf n}_\ep)\cdot \varphi\,d\sigma\right|\leq {\rm Lip}({\bf n}_\ep)\|\varphi\|^2_{L^2(\Gamma_1^\ep)^3},$$
where ${\rm Lip}({\bf n}_\ep)$ is the Lipschitz constant of the normal vector field ${\bf n}_\ep$, locally extended in a neighborhood of the surface $\{x_3=h_\ep(x')\}$ (in the sense of \cite[Section 3.4]{Boyer}). However, using definition of $h_\ep$ and condition (\ref{RelationAlpha}) on parameters $\epsilon$ and $\epsilon^\ell$, it turns out that in the general case where $\|\partial_{11}^2h\|_\infty>0$, ${\rm Lip}({\bf n}_\ep)$ is of order $\epsilon^{1-2\ell}$, hence diverging since $\epsilon^{1-\ell}$ goes to zero. As a result, we cannot use identity (\ref{remark_ref2}) to estimate in $L^2$ norms of ${\rm div}(\varphi)$ and ${\rm rot}(\varphi)$ over $\Omega^\ep$ by the $L^2$ norm of $D\varphi$, as is done in {\it step 2} of the proof of Lemma \ref{lemma_estimates} (see also \cite[Proof of Lemma 2.2]{Bayada_NewModel})  in the case of a flat boundary.
\end{remark}

\section*{Acknowledgements}
Mar\'ia dedicates this paper to her father, Julio, for all his infinite love. The authors would like to thank the anonymous referees for their nice comments that have allowed us to improve this article.

\section*{Declarations}
\paragraph{Conflict of interest} The authors declared that they have no conflict of interest.\\


\begin{thebibliography}{10}
 
\bibitem{Allaire} Allaire G.  Homogenization of the Stokes flow in a connected porous medium. Asymptot  Anal, 1989, {\bf 2}(3): 203--222

\bibitem{Anguiano_evolutivo} Anguiano M.  On the non-stationary non-Newtonian flow through a thin porous medium. ZAMM Z Angew Math Mech,  2017, {\bf  97}(8):   895--915


\bibitem{Anguiano_SG} Anguiano M, Su\'arez-Grau F J.   Nonlinear Reynolds equations for non-Newtonian thin-film fluid flows over a rough boundary. IMA J Appl Math, 2019, {\bf 84}(1):  63--95

\bibitem{Anguiano_SG_magneto} Anguiano M, Su\'arez-Grau F J.  Mathematical derivation of a Reynolds equation for magneto-micropolar fluid flows through a thin domain. Z. Angew Math Phys, 2024, {\bf 75}:  28


 

\bibitem{Bayada_Chambat_1988} Bayada G,  Chambat M.  New models in the theory of the hydrodynamic lubrication of rough surfaces. J  Tribol, 1988, {\bf 110}(3):  402--407

\bibitem{Bayada_Chambat} Bayada G,  Chambat M.   Homogenization of the Stokes system in a thin film flow with rapidly varying thickness. RAIRO Mod\'el Math Anal Num\'er, 1989, {\bf 23}(2):  205--234

\bibitem{Bayada_NewModel}   Bayada G, Benhaboucha N, Chambat M. New models in micropolar fluid and their application to lubrication. Math. Models Methods Appl. Sci., 2005, {\bf 15}(3):  343--374

\bibitem{Bayada_NewModel2}    Bayada G, Benhaboucha N, Chambat M. Wall slip induced by a micropolar fluid.  J Eng  Math, 2008, {\bf 60}:  89--100

\bibitem{Bayada_Gamouana}  Bayada g, Chambat M, Gamouana  S R.  About thin film micropolar asymptotic equations. Quart Appl Math, 2001, {\bf 59}(3): 413--439

\bibitem{Bayada_Gamouana2} Bayada G, Chambat M, Gamouana  S R. Micropolar effects in the coupling of a thin film past a porous medium.  Asymptot Anal, 2002, {\bf 30}(3--4): 187--216

 

\bibitem{BayadaLuc}   Bayada G, Lukaszewicz G. On micropolar fluids in the theory of lubrication. Rigorous derivation of an analogue of the Reynolds equation. Internat J Engrg Sci, 1996, {\bf 34}(13):  1477--1490

\bibitem{Benes}  Bene${\rm \check{s}}$ M,  Pa${\rm \check{z}}$anin I,  Radulovi\'c M, Rukavina B. Nonzero boundary conditions for the unsteady micropolar pipe flow: Well-posedness and asymptotics. Appl Math Comput, 2022 {\bf 427}:   127184


\bibitem{Benhaboucha}  Benhaboucha N,   Chambat M,  Ciuperca I.  Asymptotic behaviour of pressure and stresses in a thin film flow with a rough boundary. Quart  Appl  Math, 2005, {\bf 63}(2): 369--400
 
\bibitem{Bessonov1} Bessonov N M.  A new generalization of the Reynolds equation for a micropolar fluid and its application to bearing theory. Tribol  Int, 1994,  {\bf 27}(2): 105-108

\bibitem{Bessonov2} Bessonov N M. Boundary viscosity conception in hydrodynamical theory of lubrication, St-Petersbourg: Russian Academy of Science, Institute of the Problems of Mechanical Engineering,  1993

\bibitem{Bonn_Paz_SG2}  Bonnivard M,  Pa\v zanin I,  Su\'arez-Grau F J.  Effects of rough boundary and nonzero boundary conditions on the lubrication process with micropolar fluid. Eur J Mech B/Fluids, 2018, {\bf 72}: 501--518


\bibitem{Bonn_Paz_SG}  Bonnivard M,  Pa\v zanin I,  Su\'arez-Grau F J.   A generalized Reynolds equation for micropolar flows past a ribbed surface with nonzero boundary conditions. ESAIM: Math Model Numer Anal, 2022, {\bf 56}(4):  1255--1305


\bibitem{Boukrouche1}   Boukrouche M. A Reynolds equation rigorously derived from the micropolar Navier--Stokes system. Lecture Notes in Pure and Appl Math, 2001, {\bf 223}:  1--18

\bibitem{Boukrouche2}  Boukrouche M,  Paoli L. Asymptotic Analysis of a micropolar fluid flow in a thin domain with a free and rough boundary. SIAM J  Math  Anal, 2012, {\bf 44}(2):  1211--1256

\bibitem{Boukrouche3}  Boukrouche M,  Paoli L, Ziane F.  Unsteady micropolar fluid flow in a thin domain with Tresca fluid--solid interface law. Comput  Math  Appl, 2019 {\bf  77}(11): 2917--2932

\bibitem{Boukrouche4}  Boukrouche M,  Paoli L, Ziane F.  Micropolar fluid flow in a thick domain with multiscale oscillating roughness and friction boundary conditions. J Math  Anal  Appl, 2021, {\bf 495}(1): 124688

 

\bibitem{Boyer}  Boyer F,  Fabrie P. Mathematical Tools for the Study of the Incompressible Navier--Stokes Equations and Related Models. New York: Springer, 2013



\bibitem{Ciora}  Cioranescu D,  Damlamian A,  Griso G.  Periodic unfolding and homogenization.  C R Acad Sci Paris Ser I, 2002, {\bf 335}(1):  99--104
%
%
\bibitem{Ciora2} Cioranescu D,  Damlamian A,  Griso G.  The periodic unfolding method in homogenization. SIAM J Math Anal, 2008,  {\bf 40}(4): 1585--1620

\bibitem{Dupuy1}  Dupuy D,  Panasenko G,  Stavre R. Asymptotic methods for micropolar fluids in a tube structure. Math Models Methods Appl  Sci, 2004, {\bf 14}(5): 735--758

\bibitem{Dupuy2} Dupuy D,  Panasenko G,  Stavre R. Asymptotic solution for a micropolar flow in a curvilinear channel.  ZAMM Z  Angew  Math  Mech, 2008, {\bf 88}(10): 793--807

%

\bibitem{eringen} Eringen A.C.   Theory of micropolar fluids. J Math Mech, 1966, {\bf 16}(1): 1--18

 
%
 \bibitem{Fab2}  Fabricius J, Tsandzana A, Perez-Rafols  F,  Wall P. A Comparison of the Roughness Regimes in Hydrodynamic Lubrication, J  Tribol, 2017, {\bf 139}(5): 051702
 
\bibitem{Jacobson} Jacobson B. O.  At the boundary between lubrication and wear. In: Hutchings I M ed. First world tribology conference, London: Mech Eng Pub,  1997.  291--298

\bibitem{Leger}  L\'eger L, Hervet H,  Pit R. Friction and flow with slip at fluid-solid interfaces, Interfacial properties on the submicron scale. In:  ACS symposium series 781, ACS, Washington DC, 2001 {\bf 10}: 154--167
 
 \bibitem{Luka} Lukaszewicz G.   Micropolar fluids, theory and applications. In: Modeling and Simulation 
in Science, Engineering and Technology. Boston: Birkha$\ddot{\rm u}$ser, 1999


\bibitem{Mahabaleshwar} Mahabaleshwar U S,   Pa\v zanin I,  Radulovi\'c M, Su\'arez-Grau  F J.   Effects of small boundary perturbation on the MHD duct flow.  Theor  Appl  Mech, 2017, {\bf 44}(1):   83--101
%
\bibitem{Marusic-Paloka}  Maru\v si\'c-Paloka E, Pa\v zanin I, Maru\v si\'c   S. An effective model for the lubrication with micropolar fluid.  Mech. Res. Commun.,  2013, {\bf 52}:  69--73 

\bibitem{Marusic-Paloka2} Maru\v si\'c-Paloka E, Pa${\rm \check{z}}$anin I,   Radulovi\'c M. 
Flow of a micropolar fluid through a channel with small boundary perturbation
Z. Naturforsch. A., 2016, {\bf 71}(7)a:  607--619

\bibitem{Marusic-Paloka3} Maru\v si\'c-Paloka E, Pa${\rm \check{z}}$anin I,   Radulovi\'c M.  Justification of the higher order effective model describing the lubrication of a rotating shaft with micropolar fluid.  Symmetry, 2020, {\bf 12}(3):  334.



\bibitem{Mikelic2}  Mikeli\'c  A.  Remark on the result on homogenization in hydrodynamical lubrication by G. Bayada and M. Chambat. RAIRO Mod\'el Math Anal Num\'er, 1991, {\bf 25}(3):  363--370


\bibitem{Nakasato} Nakasato J C,  Pa\v zanin I. Homogenization of the non-isothermal, non-Newtonian fluid flow in a thin domain with oscillating boundary. Z  Angew  Math  Phys, 2023, {\bf 74}:   211


\bibitem{Pazanin}  Pa\v zanin I. Asymptotic behavior of micropolar fluid flow through a curved pipe. Acta Appl  Math, 2011, {\bf 116}:  1

\bibitem{PazaninFilomat}  Pa\v zanin I. Asymptotic analysis of the lubrication problem with nonstandard boundary conditions for microrotation. Filomat, 2016, {\bf 30}(8): 2233--2247

\bibitem{PazaninRadulovic}  Pa\v zanin I,  Radulovi\'c M. Asymptotic analysis of the nonsteady micropolar fluid flow through a curved pipe. Appl  Anal, 2020, {\bf 99}(12): 2045--2092

\bibitem{Pazanin_SG}  Pa\v zanin I, Su\'arez-Grau F J.  Analysis of the thin film flow in a rough domain filled with micropolar fluid. Comput Math Appl, 2014, {\bf 68}(12): 1915 --1932

\bibitem{Pazanin_SG_two_osci} Pa\v zanin I, Su\'arez-Grau F J.  Homogenization of the Darcy?Lapwood?Brinkman Flow in a Thin Domain with Highly Oscillating Boundaries. B Malays  Math  Sci  So, 2019, {\bf  42}:   3073 -- 3109
%
\bibitem{Pazanin_SG_thermo} Pa\v zanin I, Su\'arez-Grau F J. Roughness-induced effects on the thermomicropolar fluid flow through a thin domain. Stud  Appl  Math, 2023, {\bf 151}(2): 716--751


 

%
\bibitem{grau1} Su\'arez-Grau F J.  Asymptotic behavior of a non-Newtonian flow in a thin domain with Navier law on a rough boundary. Nonlin Anal, 2015, {\bf 117}: 99--123

\bibitem{SG1} Su\'arez-Grau F J.  Analysis of the roughness regime for micropolar fluids via homogenization. Bull Malays Math Sci Soc, 2021, {\bf 44}: 1613--1652

\bibitem{SG_porous} Su\'arez-Grau F J. Mathematical modeling of micropolar fluid flows through a thin porous medium. J  Eng Math,  2021, {\bf 126}: 7

\bibitem{SG3} Su\'arez-Grau F J. Homogenization of a micropolar fluid past a porous media with nonzero spin boundary condition. Math Meth Appl Sci,   2021, {\bf 44}(6):  4835--4857

\bibitem{Tartar}  Tartar L. Incompressible fluid flow in a porous medium convergence of the homogenization process. In: Appendix to Lecture Notes in Physics, 127. Berlin: Springer-Velag, 1980. 










%
\end{thebibliography}
\end{document}